\documentclass[10pt]{amsart}

\usepackage{hyperref}

\usepackage{amssymb}
\usepackage{amsmath}

\newcommand{\ntop}[2]{\genfrac{}{}{0pt}{1}{#1}{#2}}

\let\newpf\proof \let\proof\relax 
\newenvironment{pf}{\newpf[\proofname]}{\qed\endtrivlist}

\def\DC{{\mathrm{DC}}}

\def\bm{\begin{matrix}}
\def\em{\end{matrix}}

\newcommand{\bt}{\begin{thm}}
\newcommand{\et}{\end{thm}}

\newcommand{\bl}{\begin{lemma}}
\newcommand{\el}{\end{lemma}}

\def\ci{\cite}

\def\al{\alpha}

\newcommand{\beq}{\begin{eqnarray}}
\newcommand{\eeq}{\end{eqnarray}}

\def\bbz{{\mathbb{Z}}}

\def\bbr{{\mathbb{R}}}
\def\bbn{{\mathbb{N}}}

\def\hH{{\check H}}

\def\tn{{\underline n}}
\def\tN{{\underline N}}
\def\tphi{{\underline \phi}}
\def\tb{{\underline b}}
\def\tPhi{{\underline \Phi}}
\def\tB{{\underline B}}
\def\tPsi{{\underline \Psi}}
\def\tbeta{{\underline \beta}}

\def\EDC{{\mathrm{EDC}}}

\def\be{\begin{equation}}
\def\ee{\end{equation}}

\def\ba{{\begin{align}}}
\def\ea{{\end{align}}}

\def\d{{\underline d}}

\newcommand{\la}{\label}

\def\0{{\mathbf 0}}

\def\SL{{\mathrm {SL}}}

\def\PSL{{\mathrm {PSL}}}

\newtheorem{thm}{Theorem}[section]
\newtheorem{cor}[thm]{Corollary}

\newtheorem{lemma}[thm]{Lemma}

\theoremstyle{remark}
\newtheorem{rem}{Remark}[section]

\numberwithin{equation}{section}

\def \bn {\hfill \\ \smallskip\noindent}

\theoremstyle{definition}

\newtheorem{definition}{Definition}[section]
\def\proof{\bn {\bf Proof.} }

\def\note#1
{\marginpar

{\tiny $\leftarrow$
\par
\hfuzz=20pt \hbadness=9000 \hyphenpenalty=-100 \exhyphenpenalty=-100
\pretolerance=-1 \tolerance=9999 \doublehyphendemerits=-100000
\finalhyphendemerits=-100000 \baselineskip=6pt
#1}\hfuzz=1pt}

\newcommand{\dist}{\operatorname{dist}}

\renewcommand{\mod}{\operatorname{mod}}

\newcommand{\id}{\operatorname{id}}

\newcommand{\eps}{{\epsilon}}

\newcommand{\lb}{{\lambda}}

\newcommand{\C}{{\mathbb C}}

\newcommand{\Q}{{\mathbb Q}}
\newcommand{\R}{{\mathbb R}}

\newcommand{\Z}{{\mathbb Z}}

\def\B0{{\bold{0}}}

\catcode`\@=12

\def\Empty{}
\newcommand\oplabel[1]{
  \def\OpArg{#1} \ifx \OpArg\Empty {} \else
  	\label{#1}
  \fi}

\newcommand{\comm}[1]{}
\newcommand{\comment}[1]{}

\begin{document}

\title{Almost localization and almost reducibility}

\author{Artur Avila and Svetlana Jitomirskaya}

\address{
CNRS UMR 7599,
Laboratoire de Probabilit\'es et Mod\`eles al\'eatoires\\
Universit\'e Pierre et Marie Curie--Bo\^\i te courrier 188\\
75252--Paris Cedex 05, France
}
\curraddr{IMPA, Estrada Dona Castorina 110, Rio de Janeiro, 22460-320,
Brazil}
\email{artur@math.sunysb.edu}

\address{
University of California, Irvine, California, 92697
}
\email{szhitomi@uci.edu}
\thanks{$^\dag$This work was supported in part by NSF, grant DMS-0300974,
and BSF, grant 2002068.  This research was partially conducted during the period A.A.
served as a Clay Research Fellow.}

\begin{abstract}

We develop a quantitative version of Aubry duality and use it to obtain several sharp estimates for the dynamics of Schr\"odinger cocycles
associated to a non-perturbatively small analytic potential and
Diophantine frequency.  In particular, we
establish the full version of Eliasson's
reducibility theory in this regime (our approach
actually leads to improvements even in the perturbative regime:
we are able to show, for all energies,
``almost reducibility'' in some band of analyticity).  We also prove
$1/2$-H\"older continuity of the integrated density of states.
For the almost Mathieu operator, our results hold through
the entire regime of sub-critical coupling and imply also
the dry version of the Ten Martini Problem for the concerned parameters.

\end{abstract}

\setcounter{tocdepth}{1}

\maketitle

\section{Introduction}

This work is concerned with quasiperiodic Schr\"odinger operators $H=H_{\lambda v,\alpha,\theta}$
defined on $l^2(\Z)$
\be \label{11}
(H u)_n=u_{n+1}+u_{n-1}+\lambda v(\theta+n \alpha) u_n
\ee
where $v:\R/\Z \to \R$ is the potential, $\lambda \in \R$ is the coupling constant, $\alpha \in \R \setminus \Q$ is the
frequency and $\theta \in \R$ is the phase.  The central (and in some sense most important) example is
given by the almost Mathieu operator, when $v(x)=2\cos(2 \pi x)$.

Except where otherwise noted, below we assume the frequency $\alpha$ to
be Diophantine in the usual sense, and $v$ analytic.

An important feature of quasiperiodic operators is that the family $\{H_{\lambda v,\alpha,\theta}\}_{\lambda \in \R}$
undergoes a so called metal-insulator transition when $|\lambda|$ is changed from small to large. Loosely speaking, 
large $\lambda$ correspond to the insulator phase, with positive Lyapunov exponents and related localization type effects, while small $\lambda$ leads to a metallic phase, with zero Lyapunov exponents and good transport properties.
Here we are interested in the metallic regime, and therefore in small couplings. 

One should distinguish between two possible regimes of small $|\lambda|$
(similar considerations can be applied to the analysis of large coupling).
One is {\it perturbative}, meaning that the smallness condition on
$|\lambda|$ depends not only on the potential $v$, but also on the frequency
$\alpha$: the key resulting limitation is that the analysis at a given
coupling, however small, has to exclude a positive Lebesgue measure set of
$\alpha$. Such exclusions are inherent to the KAM-type methods that have been traditionally used in this context.
The other, stronger regime, is called {\it non-perturbative}, meaning that the
smallness condition on $|\lambda|$ only depends on the potential, leading to results that hold for
almost every $\alpha.$  Let us stress that, in certain related contexts
(multifrequency modifications of (\ref {11})),
perturbative results do fail to extend to the non-perturbative regime, see
Remark \ref {gen}.

In \cite {E}, Eliasson obtained a very precise description of operators \eqref{11}
in the case of small analytic potentials in the perturbative regime.
He proved fine estimates on {\it all}
solutions of the eigenvalue equation $H u=E u$ for $E$ in the spectrum,
concluding that most, but not all, are analytic Bloch waves, i.e.
quasiperiodic and analytic in the hull.  In his context, the problem of
existence of analytic Bloch waves can be restated (and is indeed treated)
as a dynamical systems problem, {\it reducibility} of the associated
cocycle (see \S \ref {redu def}).
His method is based on a sophisticated KAM scheme, which avoids
the limitations of earlier KAM methods (that go back to the work of
Dinaburg-Sinai \cite{ds} and that excluded parts of the
spectrum from consideration).

More recently, a less precise analysis of small analytic potentials
has been carried out through the non-perturbative regime 
 in \cite {J} and \cite {BJ1}.
One feature of those results is that most of the
analysis is concerned with the dual model $\hat H=\hat
H_{\lambda v,\alpha,\theta}$
defined on $l^2(\Z)$
\be
(\hat H \hat u)_n=\sum \lambda
\hat v_k \hat u_{n-k}+2 \cos (2 \pi \theta+n \alpha)
\hat u_n,
\ee
where $\hat v_k$ are the Fourier coefficients of $v(x)=\sum \hat v_k e^{2
\pi i k x}$.
More precisely, localization
(pure point spectrum with exponentially decaying eigenfunctions) results
for the family $\{\hat H_{\lambda v,\alpha,\theta}\}_{\theta \in \R}$
are used to obtain information on the family
$\{H_{\lambda v,\alpha,\theta}\}_{\theta \in \R}$.
At the root of this approach is
the classical Aubry duality: a Fourier-type transform corresponds
localized eigenfunctions $\hat H u=E \hat u$ to analytic Bloch waves
for the equation $H u=E u$. Duality, originally discovered in \cite{AA} has been given several
rigorous interpretations since the early 80s.
Dynamical version of Aubry duality is that 

localization for the dual model

leads to reducibility for almost every
energy \cite {puig2}. A more subtle duality statement is that pure point spectrum for the dual model allows to conclude purely absolutely continuous spectrum for a.e. $\theta$ \cite{gjls}. However, all the duality links that have been established so far were entirely algebraic, with no quantitative estimates involved. Thus  they could not lead to an
understanding of the whole spectrum, since they do not acknowledge the set
of energies where localization (and reducibility) fails. Inability to study this (zero measure yet topologically generic) set of bad energies was a major stumbling block in answering various spectral questions, including some long open conjectures. For example, a question whether there is a possibility of collapse of potential gaps, whether there is ever any singular spectrum, or what is the exact modulus of continuity of the integrated density of states all require estimates holding
for all energies.

In this paper we address this issue by developing the first {\it quantitative} version of duality, that makes it possible to  obtain fine dynamical estimates from local information on the dual model. We then show that required local properties do hold non-perturbatively for {\it all} energies in the spectrum, thus allowing us to obtain

all of Eliasson's results (and more)
in the non-perturbative setting.  As discussed above, results ``for almost
every energy'' can be obtained by considering only the localized solutions
for the dual model.  Our achievement here is in developing
the technique that can handle all energies, including those for
which localization does not hold, and there is no reducibility.  For those
energies, our approach still gives tight dynamical estimates that imply
``almost reducibility'': though coordinate changes
can not trivialize the dynamics (this would be reducibility), they come
arbitrarily close to it
(Theorem \ref{almost reducible}).  Here the notion
of closeness is rather strong, involving control in a fixed band of
analyticity (such strong control is new even in the perturbative regime).

Besides almost reducibility, the dynamical estimates yield almost
immediately a number of corollaries, including sharp estimates such as
$1/2$-H\"older continuity of the integrated density of states (Theorems \ref
{almost mathieu 1/2} and \ref {general 1/2}), and the non-collapse of
spectral gaps for
the almost Mathieu operator (Theorem \ref {dry}).  Certain other
consequences are less immediate and will be reported separately.

Our almost reducibility result obviously implies that 
the non-perturbative setting can be reduced (via coordinate change) to
the perturbative regime,\footnote {Theorem
\ref{almost reducible} is actually
much more than is needed for mere reduction to the
perturbative regime, which in itself is a less delicate result (it is
presented in the appendix as Theorem \ref{a1}).} so the known results in the perturbative
theory become automatically non-perturbative.
This is applied to conclude absolute continuity of spectral
measures for all phases (Theorems \ref {almost mathieu ac} and
\ref {general ac}) using a result of \cite{E}.  While this reduction also
gives other corollaries,
we stress that Theorems \ref{dry}, \ref{almost mathieu 1/2}, and
\ref{general 1/2} are obtained here without using this route.

\begin{rem} \label {gen}

In \cite {E}, Eliasson is actually able to deal with multifrequency
potentials.  It was shown by Bourgain \cite {B3}
that Eliasson's results, in the
multifrequency case, do not hold non-perturbatively.

\end{rem}

\begin{rem} \label {generalizations}

A further common advantage of non-perturbative approaches, which we will
not pursue here, is to extend through much weaker Diophantine
conditions than what can be covered by KAM based methods, which
usually stop working at the Brjuno condition.  The possibility of going
beyond Brjuno is very interesting for certain problems, see
for instance \cite {AJ}.

\end{rem}

\subsection{The almost Mathieu operator}

Before stating precisely the non-perturbative version of
Eliasson's result, we discuss the almost Mathieu operator. 
Results are particularly neat for this case because the
dual model of the almost Mathieu operator
is a rescaled almost Mathieu operator (with the same frequency, but inverse
coupling), thus the non-perturbative analysis extends
through the whole subcritical regime $|\lambda|<1$ (and can not be extended
further, as the almost Mathieu operator undergoes a phase
transition at $|\lambda|=1$ \cite {AA}).

In \cite {AJ}, it was proved that the spectrum of the almost Mathieu
operator is a Cantor set for any $\alpha \in \R \setminus \Q$, $\lambda \neq
0$.  This is the Ten Martini Problem of Kac-Simon.  There is a much more
difficult problem, which remains open, known as the dry version of the
Ten Martini Problem:
it asks to show that the spectrum is not only a Cantor set, but that all
gaps predicted by the Gap-Labelling Theorem \cite{blt,JM} are open
(see Section \ref{42} for a precise formulation). In contrast with the Ten Martini problem that only requires certain property to hold densely in the spectrum, this formulation requires
handling of {\it all} energies with rational rotation numbers, thus doesn't leave any room for energy exclusion. The dry Ten Martini Problem enjoyed a significant attention since its first formulation in the early 80s. An affirmative answer
was obtained for Liouville $\alpha$ \cite {CEY}, as well as for a
set of $(\lambda,\alpha)$ of positive Lebesgue measure \cite {P}.
Here we are able to deal with almost every $(\lambda,\alpha)$.

Let $\DC(\kappa,\tau)$ be the set of $\alpha \in \R$ such that
$\left |\alpha-\frac {p} {q} \right |
\geq \kappa q^{-\tau}$, $p,q \in \Z$, $q \neq 0$.  The set of
Diophantine numbers is $\DC=\cup_{\kappa>0,\tau>0} \DC(\kappa,\tau)$.

\begin{thm} \label {dry}

The dry version of the Ten Martini Problem holds for
$\alpha \in \DC$, $\lambda \neq -1,0,1$.

\end{thm}

The work of Puig mentioned above
obtains the same result in
the perturbative regime $\alpha \in \DC$, $\ln |\lambda|$ large
(depending on $\alpha$) and is based on the perturbative results of
Eliasson.

We note that our proof is not dependent on \cite{AJ}, where the issue was in handling the arithmetically critical non-Diophantine regime, while here we entirely focus on the Diophantine case.

As discussed in \cite {AJ},
in order to bridge the gap between generic and full measure one must
analyze a particularly difficult region of parameters which is still
badly understood (and which this work does not touch).   The case
$|\lambda|=1$ is also very open even in the case when $\alpha$ is the golden
mean (the problem is that one can not use localization methods).

We next move our focus to the regularity of the integrated density of
states, a recurring theme in the analysis of quasiperiodic operators.

\begin{thm} \label {almost mathieu 1/2}

Let $\alpha \in \DC$, $\lambda \neq -1,0,1$.  Then the integrated density of
states is $1/2$-H\"older continuous.

\end{thm}

This estimate is optimal in several ways. First, there are square-root
singularities at the boundaries of gaps (e.g., \cite{puig2}), so the modulus of continuity
cannot be improved. Also, it is known that for a certain non-empty set of $\alpha$ which
satisfies good Diophantine properties (but has zero Lebesgue measure)
and $\lambda=1$, the integrated density of states is not H\"older
(\cite {B2}, Remark after Corollary 8.6).
Finally, for any $\lambda \neq 0$ and generic $\alpha$, the integrated
density of states is not H\"older (this is because the Lyapunov exponent is
discontinuous at rational $\alpha$, which easily
implies that it is not H\"older for generic $\alpha$).

Goldstein-Schlag \cite {GS} had previously shown 
$(1/2-\epsilon)$-H\"older for any $\epsilon,$ all $\lambda \neq -1,0,1$ and a full
Lebesgue measure subset of Diophantine frequencies (they are also able to
consider other analytic functions, in the regime of positive
Lyapunov exponent, at the cost of a worse H\"older
constant).  Previously Bourgain \cite {B1} had obtained almost
$1/2$-H\"older
continuity in the perturbative regime, for Diophantine $\alpha$ and
$\ln |\lambda|$ large (depending on $\alpha$).  More recently
Sana Ben Hadj Amor obtained $1/2$-H\"older continuity in the perturbative
regime of Eliasson \cite{Am}.

Another consequence of our results is the following.  A well known conjecture,
dating back to the work of Aubry-Andr\'e \cite {AA}, and more recently included
in the list of problems of Simon (Problem 6 of \cite {S}),
asks to show that for
$0<|\lambda|<1$ the spectral measures of the almost Mathieu operator
are absolutely continuous. This is also tied to another, more general and far-reaching, conjecture, that is sometimes attributed to Simon, saying that for almost periodic operators, singular spectrum must be phase independent. While a.e. phase-independence of the spectrum and its a.c.,s.c. and pure point components is an almost immediate corollary of ergodicity, almost periodic operators exhibit certain phase rigidity, in that the spectrum \ci{as} and absolutely continuous 
spectrum \cite {LS} are the same for {\it all} phases. Although, this is not true for singular continuous and pure point components taken individually (there may be a dependence on the arithmetics of the phase \cite{js}), the conjecture is that combined together they will also demonstrate the phase stability. This question is only non-trivial (but highly so) in the regime of zero Lyapunov exponents. We should note, that by the time of the inclusion of the above question as Problem 6 in \cite{S}, it was already solved in case of Diophantine frequencies for a.e. phase \cite{J}, thus, in the Diophantine regime, the issue was precisely the phase rigidity. This is what we address here.

\begin{thm} \label {almost mathieu ac}

The spectral measures of the almost Mathieu operator are absolutely continuous
for $\alpha \in \DC$, $0<|\lambda|<1$, $\theta \in \R$.

\end{thm}

In the perturbative regime, $\alpha \in \DC$,
$0<|\lambda|<\lambda_0(\alpha)$, the result for all $\theta \in \R$ follows
from the work of Eliasson.  As was already mentioned,
Theorem \ref{almost mathieu ac} is obtained here from Eliasson's result
by ``non-perturbative reduction to the perturbative regime''.

\begin{rem}

Though not relevant to this work (since the estimates involved in
our reduction to the perturbative regime are already bound to the
usual Diophantine condition), we should point out that
this approach does limit possible extensions beyond the Brjuno condition
(see Remark \ref {generalizations}).  Such issues are particularly
relevant in view of recent progress towards absolutely continuous spectrum
``from the Liouville side'' \cite {ad}, which seems to break down strictly
beyond the Brjuno condition\footnote {Namely, the method of \cite {ad}
can only cover irrational numbers that are exponentially well approximated
by rationals.}, while absolutely continuous spectrum is expected to hold, as
described above, without any exceptions.\footnote {
After this work was completed, the first author has announced a complete
solution to Problem 6 of \cite {S}, which does involve the development of a
fully non-perturbative approach to absolutely continuous spectrum.}

\end{rem}

\subsection{Almost reducibility} \label {redu def}

As previously discussed, Eliasson's analysis is based on
dynamical systems considerations.  A
cocycle is defined by a pair $(\alpha,A)$ where $\alpha \in \R$ and
$A:\R/\Z \to \SL(2,\R)$ is analytic.  It is viewed as a linear skew-product
$(x,w) \mapsto (x+\alpha,A(x) \cdot w)$, $x \in \R/\Z$, $w \in \R^2$.
We say that two analytic cocycles
$(\alpha,A^{(i)})$, $i=1,2$, are analytically conjugate if there exists an
analytic map $B:\R/\Z \to \PSL(2,\R)$ such that
\be \label {conj}
A^{(2)}(x)=B(x+\alpha)A^{(1)}(x) B(x)^{-1}.
\ee
The dynamical properties of cocycles are preserved by
conjugacies.  We say that a cocycle is reducible if it is
$C^\omega$-conjugate to a cocycle of the form $(\alpha,A_*)$ where $A_*$ is
a constant matrix.
Eliasson's reducibility theory
describes the dynamics of $(\alpha,A)$ when $\alpha$ is Diophantine and $A$
is close to a constant.  He shows that such cocycles are typically
(in a measure-theoretic sense) reducible, and gives good estimates for the
non-reducible ones.  The precise closeness quantifier defines the
Eliasson's perturbative regime.  See the appendix for a summary of the
results of the theory that are relevant to this work.

Eliasson's perturbative regime is not invariant under conjugacies.  A more
intrinsic notion, almost reducibility, introduced by
Avila-Krikorian \cite {AK}\footnote {See also the first (preprint) version
of \cite {ak0}.} in the
smooth category, captures the properties of cocycles ``that behave as
cocycles in Eliasson's perturbative regime''.  The results of this work
justify introducing the corresponding definition in the analytic category
as well.

\begin{definition}

An analytic cocycle $(\alpha,A)$ is $C^\omega$-almost reducible if
the closure of its analytic conjugacy class contains a constant.

\end{definition}

If $(\alpha,A)$ is $C^\omega$-almost reducible and $\alpha$ is
Diophantine, then $(\alpha,A)$
is analytically conjugate to a cocycle in
Eliasson's perturbative regime.

The connection to Schr\"odinger operators
is clear when
\be \label {scocycle}
A(x)=S_{\lambda v,E}(x)=
\left (\bm E-\lambda v(x)&-1\\1&0 \em \right ),
\ee
since a solution $Hu=Eu$ satisfies
$A(\theta+n \alpha) \left (\bm u_n\\u_{n-1} \em \right
)=\left (\bm u_{n+1}\\u_n \em \right )$.  In other words,
the spectral properties of the family of Schr\"odinger operators
$\{H_{\lambda v,\alpha,\theta}\}_{\theta \in \R}$ are closely related to
the dynamics of the family of cocycles $\{(\alpha,S_{\lambda v,E})\}_{E \in \R}$.

\begin{thm} \label {almost reducible}

For $\alpha \in \DC$, $v:\R/\Z \to \R$ analytic,
and $0<|\lambda|<\lambda_0(v)$, the cocycles associated
with $\{H_{\lambda v,\alpha,\theta}\}_{\theta \in \R}$ are almost reducible.  If
$v(x)=2 \cos 2 \pi x$ then $\lambda_0=1$.

\end{thm}

Direct application of Eliasson's reducibility theory immediately yields
generalized version of Theorem \ref {almost
mathieu ac}:

\begin{thm} \label {general ac}

Let $\alpha,v,\lambda$ be as in Theorem \ref {almost reducible}.
Then the singular spectrum is the same
for all phases $\theta \in \R$ (and empty).

\end{thm}

We note that in the non-perturbative regime, the a.e. absolutely continuous spectrum was established in \cite{BJ1}. What we address here is the stability of singular spectrum.

\subsection{Further non-perturbative analysis}

\begin{thm} \label {general 1/2}

Let $\alpha,v,\lambda$ be as in Theorem \ref {almost reducible}.
Then the integrated density
of states is $1/2$-H\"older.

\end{thm}

The history of this question is discussed after Theorem \ref{almost
  mathieu 1/2}, and as for the almost Mathieu case, this result is
  optimal in several ways.  While we give a direct simple proof,
this theorem can also be
  derived via Theorem \ref{almost reducible} from a recent
  perturbative result of Sana Ben Hadj Amor \cite{Am}.

Our analysis also allows to study a more delicate question: H\"older
continuity of the individual spectral measures (of which the
integrated density of states is an average). We can show that for all
$\theta$ and all localized initial vectors their spectral measures are
uniformly 1/2-H\"older.  This result is once again optimal, and new even
in the perturbative regime.  This corollary
is a little more involved and also needs the introduction of some
additional theory.  It will be reported separately \cite{ajlp}.

The dry version of the Ten Martini Problem is specific to the almost
Mathieu operator, and does not hold for general analytic potentials \cite{jo}.
What can be concluded

is that for the great majority of $v$ 
(that is, excluding a set of infinite codimension),
and except for countably many $0<|\lambda|<\lambda_0(v)$, all gaps are open.

\subsection{Almost localization and quantitative Aubry duality}

Our key problem is thus to show almost reducibility of certain cocycles.

As mentioned before, there is a
classic Aubry duality link between localization and reducibility
(Theorem \ref {localreduc}), whose application is however limited since
localization (pure point spectum with exponentially decaying
eigenfunctions) in general does not hold for every $\theta$ \cite{js}. 
This of course fits with the fact that reducibility in general
does not hold for all cocycles.
(Localization for almost every $\theta$, proved in \cite {J} and \cite
{BJ1}, turns out to be enough for many results, but not for the fine ones
we are interested in here.)

It is thus natural to devise a weakened notion of localization
that could be expected to hold for every phase, and to develop ways to link it to
almost reducibility.
Here we show that this approach indeed works. Namely, we establish a quantitative version of Aubry duality
that links local exponential decay of solutions to the dual eigenvalue problem to fine dynamical estimates, thus showing that
almost reducibility of cocycles associated
$\{H_{\lambda v,\alpha,\theta}\}_{\theta \in \R}$ can be concluded from 
a property of the dual model $\{\hat H_{\lambda v,\alpha,\theta}\}$,
which we call ``almost localization'' (see Definition \ref {sle}).
Informally, almost localization gives a
precise description of the decay of generalized eigenfunctions away from
a sparse sequence of resonances, somewhat similar to what is considered in  \cite {JL},
but more precise than what is obtained in \cite {J}
and \cite {BJ1}.  Refining those results, we indeed
show in Theorem \ref {almost localized},
that almost localization holds for the dual model in the
regime relevant for Theorem \ref {almost reducible}.  Then we proceed to
show that for $\alpha \in \DC$, almost localization for the dual
model implies almost reducibility.  This link is quantitative and is 
significantly  more subtle than that
between reducibility and localization. The key dynamical estimates
needed to establish this link, Theorems \ref{triangula} and \ref{di},
coupled with Theorem \ref{almost localized}
immediately imply Theorem \ref {almost reducible} (and all its
corollaries).  They also imply directly Theorems
\ref{dry}, \ref{almost mathieu 1/2}, \ref{general 1/2}.  Those direct 
implications are not difficult but involve conceptually new arguments, 
presented in \S \ref{spectral consequences}.  Particularly, the direct 
proof of Theorem \ref{dry} is based on linking resonant rotation numbers 
to resonant phases (see Theorem \ref{link}, or Theorem \ref{almost 
reducible precise} for more detail).  We note that the proof of Theorem 
\ref{almost localized} (\S \ref{strong localization estimates}) is 
building on well developed localization methods, and borrows a number of 
ingredients from \cite{BJ1} and some shortcuts from \cite{AJ}.  Theorems 
\ref{triangula} and \ref{di} and their proofs, presented in sections
\ref{growth}, \ref{triangulariza}, and \ref{real conjugacies}, are the 
main novel technical and conceptual contributions of this paper, and are 
not based on previous work.

\begin{rem}

Previous progress in extending Eliasson's results
by coupling non-perturbative and
perturbative methods was obtained by Avila-Krikorian.  In \cite {ak0},
almost every energy results were obtained following a rather different
technique (reduction to the perturbative regime was achieved by
renormalization).  More precise (unpublished) results were later obtained
(using the connection between localization and reducibility in both
directions), including a non-perturbative proof of
reducibility under a topological condition
(which covered ``most'' of the spectrum but excluded the gap boundaries).

\end{rem}

Let us conclude with a few comments on the analysis of general operators
(\ref {11}), without restriction on the coupling constant.\footnote
{Earlier, less complete (due to reliance on
perturbative techniques), considerations (in the case of smooth potentials)
first appeared in the preprint verison of \cite {ak0}.}
There has been
much recent progress on the description of the part of the spectrum
$\Sigma_+$ corresponding to energies with a positive Lyapunov exponent (see
\cite {B2} for a lengthy account, and \cite {GS} for more recent results). 
The concept of almost reducibility allows us to determine another region of
the spectrum $\Sigma_{ar}$ which can be thoroughly analyzed (either directly
by the methods developed here, or by reduction to the perturbative regime).

There is a number of parallels between our work and previous
developments regarding the positive Lyapunov exponent regime.

For instance, the stability of positivity of the Lyapunov exponents (in other words, the 
property 
that 
the
set of energies corresponding to zero Lyapunov exponent is closed) \cite
{BJ2} is paralleled by the fact that it is also possible to conclude 
stability of 
almost reducibility from Theorem \ref
{almost reducible}, utilizing an argument of \cite{AK}.  It is precisely 
the stability of almost reduciblity that makes the concept so interesting
from a dynamical systems point of view.

It is interesting to note that almost reducibility may be after all
determined by a ``mirror'' condition to positivity of the Lyapunov exponent. 
Indeed, for energies in the spectrum, almost reducibility clearly implies
{\it strong vanishing} of the Lyaupunov exponent, in the sense that the
cocycle may grow at most subexponentially {\it in some band}.  We conjecture
that the converse also holds, that is, almost reducibility should follow
from strong vanishing of the Lyapunov exponent.

\section{Preliminaries}
For a bounded analytic function $f$
defined on a strip $\{|\Im z|<\epsilon\}$ we
let $\|f\|_\epsilon=\sup_{|\Im z|<\epsilon} |f(z)|$.  If $f$ is a bounded
continuous function on $\R$, we let $\|f\|_0=\sup_{x \in \R} |f(x)|$.
\subsection{Cocycles}

Let $\alpha \in \R \setminus \Q$, $A \in C^0(\R/\Z,\SL(2,\C))$.  We call
$(\alpha,A)$ a {\it (complex) cocycle}.
The {\it Lyapunov exponent} is given by the formula
\be
L(\alpha,A)=\lim_{n \to \infty} \frac {1} {n} \int \ln \|A_n(x)\| dx,
\ee
where $A_n$, $n \in \Z$, is defined by $(\alpha,A)^n=(n \alpha,A_n)$, so
that for $n \geq 0$,
\be
A_n(x)=A(x+(n-1)\alpha) \cdots A(x).
\ee
It turns out (since irrational rotations are uniquely ergodic), that
\be
L(\alpha,A)=\lim_{n \to \infty} \sup_{x \in \R/\Z}
\frac {1} {n} \ln \|A_n(x)\|.
\ee

\begin{rem} \label {uniform growth}

By subadditivity, for any compact set
$K \subset (\R \setminus \Q) \times C^0(\R/\Z \times \SL(2,\R))$,
for every $\delta>0$ there exists $C_{K,\delta}>0$
such that for every $k \geq 0$,
\be
\sup_{(\alpha,A) \in K}
\sup_{x \in \R/\Z} \ln \|A_k(x)\| \leq C_{K,\delta}+k
(\sup_{(\alpha,A) \in K} L(\alpha,A)+\delta).
\ee

\end{rem}

We say that $(\alpha,A)$ is {\it uniformly hyperbolic} if there exists a
continuous splitting $\C^2=E^s(x) \oplus E^u(x)$, $x \in \R/\Z$ such that
for some $C>0$, $c>0$, and for every $n \geq 0$,
$\|A_n(x) \cdot w\| \leq C e^{-c n} \|w\|$, $w \in E^s(x)$
and $\|A_{-n}(x) \cdot w\| \leq C e^{-cn} \|w\|$, $w \in E^u(x)$.  In this
case, of course $L(\alpha,A)>0$.  We say that $(\alpha,A)$ is {\it bounded}
if $\sup_{n \geq 0} \sup_{x \in \R/\Z} \|A_n(x)\|<\infty$.

Given two cocycles $(\alpha,A^{(1)})$ and $(\alpha,A^{(2)})$,
a {\it (complex) conjugacy}
between them is a continuous $B:\R/\Z \to \SL(2,\C)$ such that
(\ref {conj}) holds.
The Lyapunov exponent is clearly invariant under conjugacies.

We assume now that $(\alpha,A)$ is a {\it real} cocycle, that is,
$A \in C^0(\R/\Z,\SL(2,\R))$.  The notion of real conjugacy (between real
cocycles) is the same as before, except that we ask for $B \in
C^0(\R/\Z,\PSL(2,\R))$.  Real conjugacies still preserve the Lyapunov
exponent.  An example where it is useful to allow $B:\R/\Z \to \PSL(2,\R)$
instead of requiring $B:\R/\Z \to \SL(2,\R)$ is given by the following
well known result (see \cite {hyp} for the case of continuous
time).

\begin{thm} \label {uhyp}

Let $(\alpha,A)$ be a uniformly hyperbolic cocycle, with $\alpha$
Diophantine and $A$ analytic.  Then there exists an analytic
$B:\R/\Z \to \PSL(2,\R)$ such that $B(x+\alpha)A(x)B(x)^{-1}$ is constant.

\end{thm}

One can not always take $B:\R/\Z \to \SL(2,\R)$ in Theorem \ref
{uhyp}.\footnote{Indeed, the continuous splitting of $\R^2$ associated to a
real uniformly hyperbolic cocycle may be topologically non-trivial
(see \S 4.3 of \cite {H} for an example).}

We say that $(\alpha,A)$ is (analytically)
reducible if it is (real) conjugate to a
constant cocycle, and the conjugacy is analytic.

The fundamental group of $\SL(2,\R)$ is isomorphic to $\Z$.  Let
\be
R_\theta=\left (\bm \cos 2 \pi \theta&-\sin 2 \pi \theta\\ \sin 2 \pi
\theta&\cos 2 \pi \theta \em \right ).
\ee
Any $A:\R/\Z \to \SL(2,\R)$ is homotopic to $x \mapsto R_{nx}$ for some
$n \in \Z$ called the degree of $A$ and denoted $\deg A=n$.

Assume now that $A:\R/\Z \to \SL(2,\R)$ is homotopic to the identity.  Then
there exists $\psi:\R/\Z \times \R/\Z \to \R$ and $u:\R/\Z \times \R/\Z \to
\R^+$ such that
\be
A(x) \cdot \left (\bm \cos 2 \pi y \\ \sin 2 \pi y \em \right )=u(x,y)
\left (\bm \cos 2 \pi (y+\psi(x,y)) \\ \sin 2 \pi (y+\psi(x,y)) \em \right
).
\ee
The function $\psi$ is called a {\it lift} of $A$.  Let $\mu$ be any
probability on $\R/\Z \times \R/\Z$ which is invariant by the continuous
map $T:(x,y) \mapsto (x+\alpha,y+\psi(x,y))$, projecting over Lebesgue
measure on the first coordinate (for instance, take $\mu$ as any
accumulation point of $\frac {1} {n} \sum_{k=0}^{n-1} T_*^k \nu$ where
$\nu$ is Lebesgue measure on $\R/\Z \times \R/\Z$).  Then the number
\be
\rho(\alpha,A)=\int \psi d\mu \mod \Z
\ee
does not depend on the choices of $\psi$ and $\mu$, and is called the
{\it fibered rotation number} of
$(\alpha,A)$, see \cite {JM} and \cite {H}.

It is immediate from the definition that
\be\la{rho0}
|\rho(\alpha,B)-\theta|<C\|B-R_{\theta}\|_0.
\ee

The fibered rotation number is invariant under real conjugacies which are
homotopic to the identity.  In general, if $(\alpha,A^{(1)})$ and
$(\alpha,A^{(2)})$ are real conjugate,
$B(x+\alpha)A^{(2)}(x)B(x)^{-1}=A^{(1)}(x)$, and $B:\R/\Z \to \PSL(2,\R)$
has degree $k$ (that is, it is homotopic to $x \mapsto R_{kx/2}$) then
\be \la{rho}\rho(\alpha,A^{(1)})=\rho(\alpha,A^{(2)})+k\alpha/2.\ee

For uniformly hyperbolic cocycles there is the following
well known-result (for $\alpha$ Diophantine and $A$ analytic, it is a
consequence of Theorem \ref {uhyp}), see \S 5.17 of \cite {H}.

\begin{thm} \label {uhype}

Let $(\alpha,A)$ be a uniformly hyperbolic cocycle, with $\alpha \in \R
\setminus \Q$.  Then $2 \rho(\alpha,A) \in \alpha \Z+\Z$.

\end{thm}

\subsection{Schr\"odinger operators}

We consider now Schr\"odinger operators $\{H_{v,\alpha,\theta}\}_{\theta \in
\R}$ (we incorporate the coupling constant into $v$).  The spectrum
$\Sigma=\Sigma_{v,\alpha}$ does not depend on $\theta$,
and it is the set of $E$ such that $(\alpha,S_{v,E})$ is not uniformly
hyperbolic, with $S_{v,E}$ given by (\ref {scocycle}).

The Lyapunov exponent is defined by $L(E)=L(\alpha,S_{v,E})$.

Fixing a phase $\theta$ and $f \in l^2(\Z)$, we let $\mu=
\mu^f_{v,\alpha,\theta}$ be
the spectral measure of $H=H_{v,\alpha,\theta}$
corresponding to $f$.  It is defined so that
\be \label {resolvent}
\langle (H-E)^{-1} f,f \rangle=\int_\R \frac {1} {E'-E} d\mu(E')
\ee
holds for $E$ in the resolvent set $\C \setminus \Sigma$.  

The integrated density of states is the function $N:\R \to [0,1]$ defined by
\be
N(E)=\int_{\R/\Z} \mu^f_{v,\alpha,\theta}(-\infty,E] d\theta,
\ee
where $f \in l^2(\Z)$ is such that $\|f\|_{l^2(\Z)}=1$ (the definition is
independent of the choice of $f$).
It is a continuous non-decreasing surjective (for bounded potentials) function.
The Thouless formula relates the Lyapunov exponent to
the integrated density of states
\be
L(E)=\int_\R \ln |E'-E| dN(E').
\ee
There is also a relation to the fibered rotation number
\be
N(E)=1-2 \rho(\alpha,S_{v,E})
\ee
where $\rho(\alpha,S_{v,E}) \in [0,1/2]$.

\subsection{The dual model}

It turns out that the spectrum $\hat \Sigma$ of $\hat H=\hat
H_{v,\alpha,\theta}$ coincides with the spectrum $\Sigma$ of
$H=H_{v,\alpha,\theta},$ a manifestation of Aubry duality (e.g., \cite{gjls}).
The spectral measures
$\hat \mu=\hat \mu^f_{v,\alpha,\theta}$ can be defined using the analogous
formula to (\ref {resolvent}).  The integrated density of states $\hat
N(E)=\int \hat \mu_\theta(-\infty,E] d\theta$ then coincides with
  $N(E)$.
Berezanskii's theorem \ci{ber,sim} gives in this context:

\begin{thm} \label {generalized eigenfunctions}

For any $v,\alpha,\theta,f$, and $\hat \mu^f_{v,\alpha,\theta}$-almost every
$E$, there exists a non-zero
solution $\hat H \hat u=E \hat u$ with $|\hat u_k| \leq 1+|k|$.

\end{thm}

Since $\Sigma=\hat \Sigma$ is the union of the supports of the spectral
measures, it follows:

\begin{thm} \label {linear growth}

For any $v,\alpha,\theta$, there exists a dense set of $E \in
\Sigma=\Sigma_{v,\alpha}$
such that there exists a non-zero solution $\hat H \hat u=E \hat u$ with
$|\hat u_k| \leq 1+|k|$.

\end{thm}

\subsection{ Localization and reducibility: dynamical formulation of Aubry duality. }

We describe the 
connection between localization and reducibility we mentioned in the
introduction. The Theorem below is essentially a careful dynamical formulation of the
classical Aubry duality. It appears in a similar form in \cite{puig2}.

\begin{thm} \label {localreduc}

Let $\alpha \in \R \setminus \Q$ and let
$v:\R/\Z \to \R$ be analytic.  Let
$\theta \in \R$, $E \in \R$ be such that there
exists a non-zero exponentially decaying solution
$\hat H \hat u=E \hat u$, and let $A=S_{v,E}$.
\begin{enumerate}
\item If $2 \theta \notin \alpha \Z+\Z$, there exists $B:\R/\Z \to
\SL(2,\R)$ analytic such that $B(x+\alpha) A(x) B(x)^{-1}=R_{\pm \theta}$
(so $(\alpha,A)$ is reducible). 
In particular $\|A_n(x)\|=O(1)$, $x \in \R/\Z$.
\item If $2\theta-k\alpha \in \Z$ for some $k \in \Z$, there exists
$B:\R/\Z \to \PSL(2,\R)$ and $\kappa:\R/\Z \to \R$ analytic such that
$B(x+\alpha) A(x) B(x)^{-1}=\left (\bm \pm 1&\kappa(x)\\0& \pm 1 \em
\right )$. 
In particular $\|A_n(x)\|=O(n)$, $x \in \R/\Z$.
If moreover $\alpha \in \DC$,
then $\kappa$ can be chosen to be a constant (and $(\alpha,A)$ is
reducible).
\end{enumerate}
In either case $\rho(\alpha,S_{v,E})=\pm\theta+m \alpha/2,$ for some $m\in\Z.$
\end{thm}

\begin{pf}

Let $u(x)=\sum \hat u_k e^{2 \pi i k x}$ and $U(x)=\left
(\bm e^{2 \pi i \theta} u(x)\\u(x-\alpha) \em \right )$.
Then $A(x) \cdot U(x)=e^{2
\pi i \theta} U(x+\alpha)$.  Let $\tilde B(x)$ be the matrix with columns
$U(x)$ and $\overline U(x)$.  By minimality of $x \mapsto x+\alpha$,
$\det \tilde B$ is a constant.

(Case A)\,
If $\det \tilde B \neq 0$, we have $\tilde B(x+\alpha)^{-1} A(x) \tilde
B(x)=\left (\bm e^{2 \pi i \theta}&0\\0&e^{-2 \pi i \theta} \em \right )$.
It is easy to see that $\det \tilde B=\pm c i$ for some $c>0$, we then take
$B^{-1}=\frac {1} {(2c)^{1/2}}
\tilde B \left (\bm 1&\pm i\\1&\mp i \em \right )$.

(Case B)\, If $\det \tilde B=0$ then
$U(x)=\psi(x) W(x)$ with $W(x)$ a real vector defined up to sign and
$|\psi(x)|=1$.  By
minimality of $x \mapsto x+\alpha$, $W(x) \neq 0$ for every $x \in \R/\Z$.
The matrix $B(x)^{-1} \in \PSL(2,\R)$
with columns $W$ and $\frac {1} {\|W(x)\|^2} R_{1/4} W(x)$ is
thus well defined.
Then $B(x+\alpha)A(x)B(x)^{-1}=
\left (\bm d(x)&\kappa(x)\\0& d(x)^{-1}\em \right )$, where $d(x)=\frac
{\psi(x+\alpha)} {\psi(x)} e^{2 \pi i \theta}$.  Since
 $|d(x)|=1$ and $d(x)$ is real, $d(x)=\pm 1$.  If $\alpha \in \DC$,
we can also further conjugate $A$ to a constant parabolic (or identity)
matrix by solving (using, say, Fourier series) the cohomological equation
$\pm \phi(x+\alpha) \mp \phi(x)=\kappa(x)-\int_0^1
\kappa(x) dx$ with $\int_0^1 \phi(x) dx=0$ in $\R/\Z$.  Letting $B'(x)=
\left (\bm 1&-\phi(x)\\0&1 \em \right ) B(x)$, we get
$B'(x+\alpha) A(x) B'(x)^{-1}=
\left (\bm \pm 1&\int_0^1 \kappa(x) dx\\0&\pm 1 \em \right )$.

Assume that $2 \theta \notin \alpha \Z+\Z$.  Then we can not be in
case (B): indeed $e^{2 \pi i \theta}=\pm \frac {\psi(x)} {\psi(x+\alpha)}$
implies (using Fourier series) that
$\psi(x)=e^{-\pi i k x}$ (notice that $\psi$ is well defined only in
$\R/2\Z$) and $e^{2 \pi i \theta}=\pm e^{\pi i k \alpha}$, that is $2
\theta=k\alpha \mod \Z$.  Thus we are in case (A), and the first statement
follows immediately.

Assume now that $2\theta-k\alpha \in \Z$.  If we are in case (B), the second
statement follows immediately.  Otherwise we are in case (A), and we have
$B(x+\alpha) A(x) B(x)^{-1}=R_{\pm \theta}$ for some $B:\R/\Z \to
\SL(2,\R)$.  We then set $B'(x)=R_{\mp \frac {kx} {2}} B(x)$, and we get
$B'(x+\alpha)A(x)B'(x)^{-1}=\pm \id$.  So the second statement still follows
in this case (with $\kappa(x)=0$).

The statement about the value of $\rho$ follows immediately from (\ref{rho}).
\end{pf}

If $2 \theta \in \alpha \Z+\Z$, we will say that $\theta$ is {\it rational}
(with respect to $\alpha$).

\begin{rem} \label {qua}

It is clear from the above proof that if the Fourier transform of $\hat u$
is analytic in a strip $|\Im z|<\epsilon$ then the matrix $B$ given in item
(1) is analytic in the same strip.  As for item (2), it is still possible to
define a conjugating matrix with a definite complex extension, but one must
be more careful in its definition.

\end{rem}

\subsection{Rational approximations}

Let $q_n$ be the denominators of the
approximants of $\alpha$.  We recall the basic properties:
\be \label {b1}
\|q_n \alpha\|_{\R/\Z}=\inf_{1 \leq k \leq q_{n+1}-1}
\|k\alpha\|_{\R/\Z},
\ee
\be \label {b2}
1 \geq q_{n+1} \|q_n \alpha\|_{\R/\Z} \geq 1/2.
\ee

One aspect of the ``good distribution'' of orbits
$\{x+j\alpha\}_{j=0}^{q_n-1}$ which will find repeated use in this work
is the following estimate.

\begin{lemma}[Lemma 9.7 of \cite {AJ}] \label {9.7}

Let $\alpha \in \R \setminus \Q$,
$x \in \R$ and let $0 \leq l_0 \leq q_n-1$ be such that $|\sin \pi
(x+l_0 \alpha)|$ is minimal.  Then for some absolute constant $C>0$,
\be \label {irrat-1}
-C \ln q_n \leq \sum_{\ntop {l=0} {l \neq l_0}}^{q_n-1}
\ln |\sin \pi (x+l \alpha)|+(q_n-1) \ln 2 \leq C \ln q_n.
\ee

\end{lemma}

\section{General setup and statements of the main (localization and
dynamical) estimates}

\subsection{Almost localization for every $\theta$}\label{3.1}

Let $\alpha \in \R$, $\theta \in \R$, $\epsilon_0>0$.
We say that $k$ is an $\epsilon_0$-resonance if $\|2
\theta-k\alpha\|_{\R/\Z} \leq e^{-|k|\epsilon_0}$ and $\|2
\theta-k\alpha\|_{\R/\Z}=\min_{|j| \leq |k|} \|2 \theta-j\alpha\|_{\R/\Z}$.

\begin{rem}

In particular, there always exists at least one resonance, $0$.  If $\alpha
\in \DC(\kappa,\tau)$, $\|2
\theta-k\alpha\|_{\R/\Z} \leq e^{-|k|\epsilon_0}$ implies $\|2
\theta-k\alpha\|_{\R/\Z}=\min_{|j| \leq |k|} \|2 \theta-j\alpha\|_{\R/\Z}$
for $k>C(\kappa,\tau)$.

\end{rem}

We order the
$\epsilon_0$-resonances $0=n_0<|n_1| \leq |n_2| \leq ...$.  We say that $\theta$
is $\epsilon_0$-resonant if the set of resonances is infinite.
If $\theta$ is non-resonant, with the set of resonances
$\{n_0,\ldots,n_j\}$ we formally set $n_{j+1}=\infty.$  The
Diophantine condition immediately implies exponential repulsion of
resonances:

\begin{lemma} \label {following resonance}

If $\alpha \in \DC(\kappa,\tau)$, then
$|n_{j+1}| \geq c \|2 \theta-n_j \alpha\|_{\R/\Z}^{-c} \geq
c e^{c \epsilon_0 |n_j|}$, where $c=c(\kappa,\tau,\epsilon_0)>0$.

\end{lemma}

\begin{rem} \label {especially}

In case $\|2\theta-\ell_0 \alpha\|_{\R/\Z}=0,$ (so that $\theta$ is rational
with respect to $\alpha$) we have an
especially strong resonance at $\ell_0$.  In particular, $\theta$ is
non-resonant, as there will be no resonances $n$ with $|n|>|\ell_0|$.

\end{rem}

A simple Borel-Cantelli argument shows that the set of non-resonant $\theta$
has full Lebesgue measure.

\begin{definition} \label {sle}

We say that the family $\{\hat H_{v,\alpha,\theta}\}_{\theta \in \R}$
exhibits almost localization  if there exist $C_0>0$, $C_1>0$,
$\epsilon_0>0$, $\epsilon_1>0$,
such that for every solution $\hat u$ of $\hat H_{v,\alpha,\theta} \hat u=E \hat u$
satisfying $\hat u_0=1$ and $|\hat u_k| \leq 1+|k|$, and for every $C_0
(1+|n_j|)<k<C_0^{-1} |n_{j+1}|$ we have $|\hat u_k| \leq C_1 e^{-\epsilon_1 k}$
(where the $n_j$ are the $\epsilon_0$-resonances of $\theta$).

\end{definition}

\begin{rem} \label {localization non-resonant}

It is clear from Theorem \ref {generalized eigenfunctions} that almost
localization implies localization for non-resonant $\theta$
(slowly growing generalized eigenfunctions can be always shifted and
normalized to match the definition).

\end{rem}

\begin{thm} \label {almost localized}

If $v:\R/\Z \to \R$ is analytic and $|\lambda|<\lambda_0(v)$\footnote
{ The proof actually gives a quantitative bound on the
dependence of $\lambda_0$ on the analytic extension of $v$.  More precisely,
what is needed is that $\sup_{|\Im x|<\epsilon} |\lambda_0 v(x)| \leq c_0
\epsilon^{k_0}$ for some $0<\epsilon<1$, where $c_0$ and $k_0$ are absolute
constants, see \eqref{cond1},\eqref{cond2},\eqref{518}.} then $\{\hat
H_{\lambda v,\alpha,\theta}\}_{\theta \in \R}$ is almost localized for every $\alpha \in \DC$.
For $v(x)=2 \cos 2 \pi x$, we have $\lambda_0=1$.

\end{thm}

This theorem will be proved in \S \ref {strong localization estimates}.

\subsection{Bounded eigenfunctions for every energy}

The next result allows to pass from ``every $\theta$'' statements to ``every
$E$'' statements.

\begin{thm} \label {boundedeigenfunctions}

If $E \in \Sigma$
then there exists $\theta \in \R$ and a bounded solution of
$\hat H_{v,\alpha,\theta} \hat u=E \hat u$ with $\hat u_0=1$ and
$|\hat u_n| \leq 1$.

\end{thm}

\begin{pf}

It is enough to show this for a dense set of $E$ in the spectrum.
By Theorem \ref {linear growth}, we may assume that
there is a generalized eigenfunction $\hat u'$
with subexponential growth for some
phase.  Fix $N>0$, $\epsilon>0$.
Let $k_i$ be a sequence such that $|\hat u'_{k_i}|=\max_{|j| \leq (i+1) N}
|\hat u'_j|$.  It follows that there
exists some (and indeed infinitely many)
$i$ with $|\hat u'_{k_{i+1}}| \leq (1+\epsilon)
|\hat u'_{k_i}|$.  By shifting the
phase and rescaling, we obtain, for every $\epsilon>0$, $N>0$, some phase
$\theta^{N,\epsilon}$
such that there exists an eigenfunction with $\hat u^{N,\epsilon}_0=1$
and $|\hat u^{N,\epsilon}_n| \leq 1+\epsilon$, $|n| \leq N$.
Passing to the limit as $N \to \infty$ and $\epsilon \to 0$

we get
the desired eigenfunction.\footnote {An alternative proof is the following.
Let $p/q$ be close to $\alpha$, and let $\theta' \in \R$.  The spectrum of
$\hat H'=\hat H_{v,p/q,\theta'}$
is close to the
spectrum of $\hat H$ in the Hausdorff topology.  Let $E'$ be close to
$E$ and in the spectrum of $\hat H'$.  Then there is a
non-zero periodic solution to $\hat H' \hat u'=E' \hat u'$.  Changing
$\theta'$ to $\theta'+k p/q$,
we may assume that $\hat u'_0=1$ and $|\hat u'_n| \leq
1$.  Taking the limit $p/q \to \alpha$, $E' \to E$, a limit $\theta$ of
$\theta'$, and a pointwise limit $\hat u$ of $\hat u'$,
we get the statement.}
\end{pf}

\begin{rem}

The result can be generalized to some classes of continuous ergodic
operators.  For dynamical systems, corresponding results appear in the work
of Ma\~n\'e.  Namely, if a cocycle is not uniformly hyperbolic then there
exists a vector that is never expanded, either in the future or in the past.

\end{rem}

\subsection{Main dynamical estimates: quantitative Aubry duality} \label {general setup}

We fix $\alpha,v,\lambda$ as in Theorem \ref {almost localized}. 
For every $E \in \Sigma_{\lambda v,\alpha}$, let $\theta=\theta(E)$ be given
by Theorem \ref {boundedeigenfunctions}.\footnote {\label {foot}
Notice that $\theta(E)$
is not necessarily uniquely defined.  There could be uncountably many
resonant $\theta$, but at most finitely many
non-resonant $\theta$, corresponding to the same $E$.  This does not concern
our arguments, we just fix some $\theta$.}
Let $\{n_j\}$ be the set of resonances of $\theta(E)$.  Let $A=S_{\lambda
v,E}$.  In what follows,
$C$ is a large constant and $c$ is a small constant, which are allowed to
depend on $v,\lambda,\alpha$, but not on $E$ or $\theta$.

\begin{thm} \label {triangula}

Fix some $n=|n_j|+1<\infty$ and let $N=|n_{j+1}|$.
Then there exists $\Phi:\R/\Z \to \SL(2,\C)$
analytic with $\|\Phi\|_{c n^{-C}} \leq C n^C$ such that
\be
\Phi(x+\alpha)A(x)\Phi(x)^{-1}=\left (\bm e^{2 \pi i \theta} & 0\\0 & e^{-2
\pi i \theta} \em \right )+\left (\bm q_1(x)&q(x)\\
q_3(x)&q_4(x) \em \right ),\ee
with
\be \label {q1}
\|q_1\|_{c n^{-C}},\|q_3\|_{c n^{-C}},\|q_4\|_{c n^{-C}}
\leq C e^{-c N}
\ee
and
\be \label {q}
\|q\|_{c n^{-C}} \leq C e^{-cn (\ln (1+n))^{-C}}.
\ee

\end{thm}

\begin{thm} \label {di}

Fix some $n=|n_j|+1<\infty$ and let $N=|n_{j+1}|$.  Let
$L^{-1}=\|2 \theta-n_j \alpha\|_{\R/\Z}$, and assume that
$0<L^{-1}<c$.\footnote {It is likely
that the result holds assuming only $L^{-1}>0$.} 
Then there exists $W:\R/\Z \to \SL(2,\R)$ analytic such that $|\deg W|
\leq C n$, $\|W\|_c \leq C L^C$ and $\|W(x+\alpha) A(x)
W(x)^{-1}-R_{\mp \theta}\|_c \leq C e^{-c N}$.

\end{thm}

\begin{rem}

For $N=\infty,$ Theorem \ref {di} gives a
quantitative version of the first statement of Theorem \ref{localreduc}.

\end{rem}

We will prove Theorem \ref {triangula} in \S \ref {triangulariza} and
Theorem \ref {di} in \S \ref {real conjugacies}.  All spectral results are
consequences of those theorems combined with Theorem \ref{almost localized}.

\subsection{Outline of the rest of the paper}

Almost reducibility (Theorem \ref {almost reducible}),
and the direct proof of Theorem \ref {dry} are immediate consequences of
Theorem \ref {di}.  Estimates related to modulus of continuity,
including Theorem \ref {general 1/2}, are obtained from Theorem \ref
{triangula}.  All are more or less immediate.  We will discuss those
consequences in \S \ref {spectral consequences}.

The technical core of the paper is formed by the proofs of
Theorems \ref {almost localized}, \ref {triangula} and \ref
{di}.  The dynamical estimates build on preliminary estimates and ideas
developed in \S \ref {growth}, but are otherwise independent.
Localization (\S \ref{strong localization estimates}) and quantitative duality (\S\S\ref{growth},\ref{triangulariza}, and \ref{real conjugacies})  are independent. Section \ref{strong localization estimates} uses the machinery developed in \cite{BJ1} and some shortcuts from \cite{AJ}. 
The techniques and ideas developed in \S\S \ref {growth},\ref{triangulariza}, and \ref{real conjugacies} are new and do not use any ideas/methods from the existing literature.

\section{Easy spectral consequences of the main dynamical estimates} \label
{spectral consequences}

\subsection{Almost reducibility}

We will show the following precise version of
Theorem \ref {almost reducible}.

\begin{thm} \label {almost reducible precise}

Assume $\alpha \in \DC, 0<|\lambda|<\lambda_0(v).$ There exists $c>0$ (depending on
$\lambda,v,\alpha$) with the following
property.  Let $A=S_{\lambda v,E}$.

\begin{enumerate}

\item If $\rho(\alpha,A)$ is $c$-resonant, there exists a sequence
$B^{(n)}:\R/\Z \to \SL(2,\R)$ such that $B^{(n)}(x+\alpha) A(x) B^{(n)}(x)^{-1}$ converges
to a constant rotation uniformly in $\{|\Im x|<c\}$.

\item If $\rho(\alpha,A)$ is not $c$-resonant and $2 \rho(\alpha,A) \notin
\alpha \Z+\Z$ then there exists $B:\R/\Z \to \SL(2,\R)$, analytically
extending to $\{|\Im z|<c\}$,
such that $B(x+\alpha) A(x) B(x)^{-1}$ is a constant rotation.

\item If $2 \rho(\alpha,A) \in
\alpha \Z+\Z$ then there exists $B:\R/\Z \to \PSL(2,\R)$ analytic
such that $B(x+\alpha) A(x) B(x)^{-1}$ is a constant.

\end{enumerate}

\end{thm}

\begin{pf}

If $E \notin \Sigma$, by Theorem \ref {uhype},
$2 \rho(\alpha,A) \in \alpha \Z+\Z$,
and by Theorem \ref {uhyp}, $(\alpha,A)$ is reducible.

Let $E \in \Sigma$.  If $\theta$ is not $\epsilon_0$-resonant,
by Theorem \ref{almost localized} there is an
exponentially decaying eigenfunction.  Theorem \ref {localreduc} thus
applies, and the result holds in all cases (using Remark \ref {qua}).

Assume that $\theta$ is $\epsilon_0$-resonant (and thus $2\theta \notin
\alpha \Z+\Z$).  Applying Theorem \ref {di}, we get
for every $j$ a matrix $B^{(j)}:\R/\Z \to \SL(2,\R)$ such that
$\Phi(x)=B^{(j)}(x+\alpha)A(x)B^{(j)}(x)^{-1}$ satisfies
$\|\Phi(x)-R_{\mp \theta}\|_c \leq Ce^{-c N}$ and $|\deg B^{(j)}| \leq C
(|n_j|+1)$ (where $N=|n_{j+1}|$).
To conclude, let us show that $\rho$ is $c$-resonant.

By (\ref{rho}), $|\rho(\alpha,A)\pm \theta+\deg B^{(j)} \alpha|\leq Ce^{-cN}$.
Using (\ref{rho0}) and Lemma \ref {following resonance} we get for large $j$
\be \label {minorat}
\|2 \rho(\alpha,A)\pm(n_j \pm 2 \deg B^{(j)}) \alpha\|_{\R/\Z}\geq
\|\mp 2 \theta \pm n_j \alpha\|_{\R/\Z}-C e^{-c N} \geq 
c N^{-C}-C e^{-c N}>0,
\ee
\begin{align} \label {majorat}
\|2 \rho(\alpha,A)\pm(n_j \pm 2 \deg B^{(j)}) \alpha\|_{\R/\Z} &\leq
\|\mp 2 \theta \pm n_j \alpha\|_{\R/\Z} +C e^{-c N}\\
\nonumber
&\leq e^{-\epsilon_0 |n_j|}+C e^{-c N}
\leq e^{-c|n_j|} e^{-c (|n_j \pm 2 \deg B^{(j)}|)}.
\end{align}
For $j$ large, (\ref {majorat}) and the Diophantine condition
imply that $\rho(\alpha,A)$ has a $c$-resonance at
$\mp (n_j \pm 2 \deg B^{(j)})$.
If the set of $c$-resonances for $\rho(\alpha,A)$
would not be infinite, then
we would have $\|2 \rho(\alpha,A)\pm(n_j \pm 2 \deg B^{(j)}) \alpha\|_{\R/\Z}=0$
for large $j$.  But this contradicts (\ref {minorat}).
\end{pf}

\begin{rem}

This result is new even in the perturbative regime.  Though Eliasson's
scheme provides a sequence of approximate conjugacies to normal forms, for
badly behaved $\rho$ there is only control in a
shrinking sequence of bands, as typical for KAM analysis. In fact it does not imply our strong
definition of almost reducibility of cocycles except when it actually
implies reducibility.

\end{rem}

\begin{rem}

The analytic conjugacy given in the
third statement for $E \in \Sigma$ can be shown to extend also to a definite
strip.  For $E \notin \Sigma$, we do not get any estimates on the
analytic extension.

\end{rem}

We also state separately the following statement, already obtained as
a part of the proof of Theorem \ref{almost
reducible precise}. By Theorem \ref{localreduc}, if $\theta(E)$ (as
specified in Theorem \ref{boundedeigenfunctions}, see also the
footnote \ref {foot} in Section \ref{general setup}) is not
$\epsilon_0$-resonant, then $\rho(E)=\pm\theta+k\alpha/2$ for some
$k.$
There are reasons to believe that this does not hold for {\bf all} $E$ or $\theta.$\footnote{
Otherwise the set of $E$'s corresponding to a given $\theta$ would be
countable even when $\hat H_\theta$ has singular continuous spectrum.
If we relaxed slightly the definition of $\theta(E)$ to require
$|\hat u_0|>1-\epsilon$ instead of $|\hat u_0|=1$ we could have argued that
this is a
  contradiction since it can be shown in certain cases that all
  generalized eigenfunctions are bounded, and therefore for all
  energies $E$ in the support of singular-continuous measure for a given
  $\theta$
one can find $\theta(E)$ of the form $\theta+k\alpha$ for some $k.$ Such relaxed definition of $\theta(E)$ could also be used in our proofs, with only small 
changes.}

Instead, for resonant $\theta$ we have the following statement.

\begin{thm}\label{link}
Fix $E\in \Sigma$ and some $\theta(E).$ Under the conditions of
Theorem \ref{almost reducible precise} and with the same $c$ we have
that
if $\theta$ is $\epsilon_0$-resonant, then $\rho$ is $c$-resonant. 
\end{thm}

\subsection{Open gaps for the almost Mathieu operator (direct proof)}\label{42}

By Theorem \ref {uhype}, in the closure of a
component of $\R \setminus \Sigma$ we must
have $N(E)=1-2 \rho \in \alpha \Z+\Z$.
The dry Ten Martini Problem is the conjecture that 
the converse holds in the case of the almost Mathieu operator.  Since the
integrated density of states is the same for $H$ and $\hat H$
it doesn't matter whether
to prove such statement for $\lambda$ or for $\lambda^{-1}$.
Theorem \ref {dry} is thus a consequence of the following.

\begin{thm}

Let $v(x)=2 \cos 2 \pi x$, and let $0<|\lambda|<1$,
$\alpha$ Diophantine.  If $E \in \Sigma$
is such that $N(E) \in \alpha \Z+\Z$ then $E$ belongs to
the boundary of a component of $\R \setminus \Sigma$.

\end{thm}

\begin{pf}

By \cite {P}, it is enough to show that $(\alpha,S_{\lambda v,E})$
is reducible.  By Theorem \ref {almost reducible precise}, if this was not
the case then $\rho(\alpha,S_{\lambda
v,E})$ would be $c$-resonant.  By the Diophantine condition and
Remark \ref {especially}, $2 \rho \notin \alpha \Z+\Z$.  Since
$N=1-2\rho$, this is a contradiction.

\end{pf}

\subsection{Precise bounds on growth, complex perturbations}

Conventions below are as in section \ref {general setup}.

\begin{thm} \label {epsilon14}

Let $n=|n_j|+1<\infty$ and let $N=|n_{j+1}|.$
There exists $W=W_{(\epsilon)}:\R/\Z \to \SL(2,\C)$
analytic such that, letting
$Z(x)=Z_{(\epsilon)}(x)=W_{(\epsilon)}(x+\alpha)A(x)W^{-1}_{(\epsilon)}(x)$,
we have
\be \label {epsilonC}
\|W\|_{c n^{-C}} \leq C \epsilon^{-1/4} \quad \text {and} \quad
\|Z\|_{c n^{-C}} \leq 1+C e^{-c n (\ln (n+1))^{-C}}
\epsilon^{1/2}, \quad \text {for} \quad
C e^{-c N} \leq \epsilon \leq c n^{-C}.
\ee

\end{thm}

\begin{pf}

Let $\Phi$ be given by Theorem \ref {triangula}.  Let $D=\left (\bm
d&0\\0&d^{-1} \em \right )$ where
$d=\|\Phi\|_{c n^{-C}} \epsilon^{1/4}$.  Let $W(x)=D \Phi(x)$.  If
$\epsilon \leq c n^{-C}$ we have $\|W\|_{c n^{-C}} \leq C \epsilon^{-1/4}$.
Since 
\be \label{d}
D\left (\bm
a_1&a_2\\a_3&a_4 \em \right )D^{-1}=
\left (\bm
a_1&d^{2}a_2\\d^{-2}a_3&a_4 \em \right )\ee we get
\be
Z(x)=W(x+\alpha)A(x)W(x)^{-1}=\left (\bm e^{2 \pi i \theta} & 0\\0 & e^{-2
\pi i \theta} \em \right )+\left (\bm z_1(x)&z_2(x)\\
z_3(x)&z_4(x) \em \right ),
\ee
with $\|z_1\|_{c n^{-C}},\|z_3\|_{c n^{-C}},\|z_4\|_{c n^{-C}}
\leq C \epsilon^{-1/2} e^{-c N}$ and $\|z_2\|_{c n^{-C}} \leq C
\epsilon^{1/2} e^{-cn (\ln (n+1))^{-C}}$.
If $\epsilon \geq C e^{-c N}$ then
$\|Z\|_{c n^{-C}} \leq 1+C e^{-cn(\ln (n+1))^{-C}} \epsilon^{1/2}$.
\end{pf}

The following gives a direct proof of a perturbative estimate of \cite {E} (see Theorem \ref{eli} in the Appendix).

\begin{cor} \la{linear}

For every $s \geq 0$ we have $\|A_s\|_0 \leq C (1+s)$.  Moreover, if
$\theta$ is non-rational then $\|A_s\|_0=o(1+s)$.

\end{cor}

\begin{pf}

By Lemma \ref {following resonance},
for every $s \geq C$, $\epsilon=1/s^2$ is in the range specified in (\ref
{epsilonC}) for some choice of $n=|n_j|+1$.  Then
\be
\ln \|A_s\|_0 \leq 2 \ln \|W\|_0+s \ln \|Z\|_0 \leq C-\frac {1} {2}
\ln \epsilon+s C e^{-c n(\ln (n+1))^{-C}} \epsilon^{1/2} \leq C+\ln (1+s).
\ee
This gives $\|A_s\|_0 \leq C(1+s)$.

If $\theta$ is non-rational and non-resonant then by Theorems \ref{almost
localized}, \ref{localreduc}, $\|A_s\|=O(1)$.  If
$\theta$ is resonant then we have $n \to \infty$ as $s \to \infty$, and we get
the estimate
\be
\sup_{0 \leq j \leq c e^{c n (\ln n)^{-C}} s} \ln \|A_j\| \leq C+\ln (1+s),
\ee
which implies the second statement.
\end{pf}

\begin{cor} \label {perturbati}

If $B:\R/\Z \to \SL(2,\C)$ is continuous then
$L(\alpha,B) \leq C \|B-A\|_0^{1/2}$.

\end{cor}

\begin{pf}

It is enough to consider the case when $\epsilon=\|B-A\|_0$ is sufficiently small.  Then
$\epsilon$ is in the range specified by (\ref {epsilonC}) for some
$n=|n_j|+1$.  Let $\tilde B(x)=W(x+\alpha) B(x) W(x)^{-1}$.  Then
$\|\tilde B\|_0 \leq \|Z\|_0+\|W\|_0^2 \|B-A\|_0 \leq 1+C \epsilon^{1/2}$.
Then $L(\alpha,B)=L(\alpha,\tilde B) \leq \ln \|\tilde B\|_0 \leq C
\epsilon^{1/2}$.
\end{pf}

As before, the estimate is improved when $\theta$ is not rational.

\subsubsection{The integrated density of states:
direct proof of Theorem \ref {general 1/2}}

Theorem \ref {general 1/2} follows easily from Corollary \ref {perturbati}.
Indeed by the Thouless formula
$L(E)=\int \ln |E-E'| dN(E')$.  Thus $L(E+i\epsilon) \geq
L(E+i\epsilon)-L(E)=\frac {1} {2}
\int\ln(1+\frac{\epsilon^2}{(E-E')^2})dN(E')\geq c
(N(E+\epsilon)-N(E-\epsilon))$ for every $\epsilon>0$.  By Corollary \ref
{perturbati}, $L(E+i\epsilon) \leq C \epsilon^{1/2}$, $E \in \Sigma$. 
Thus $N(E+\epsilon)-N(E-\epsilon) \geq C c^{-1} \epsilon^{1/2}$ for every
$0<\epsilon<1$, $E \in \Sigma$.  Since $N$ is locally constant in the
complement of $\Sigma$, this means precisely that $N$ is $1/2$-H\"older.

\section{Almost localization: proof of Theorem \ref {almost localized}}
\label {strong localization estimates}

Although it will not be needed for the rest of the paper we will consider a weaker Diophantine condition on $\alpha$.
For $\nu>0$, $\xi>0$, let $\EDC(\nu,\xi)$ be the set of $\alpha$ such that
\be
|q\alpha-p| \geq \nu e^{-\xi q}, \quad p \in \Z, q \in \Z \setminus
\{0\}.
\ee
Clearly for all $\kappa,\tau,\xi>0$ there exists $\nu>0$ such that
$\DC(\kappa,\tau) \subset \EDC(\nu,\xi)$.  If $\alpha \in \EDC(\nu,\xi)$
then, by (\ref{b2}), $\ln q_{n+1} \leq \xi q_n-\ln \nu$.

We will prove the
following precise version of Theorem \ref {almost localized}.
\bt \la{SLE}
%SLE for any $C>1$. 
There exists $\lambda_0(v)>0$ such that if
 $\,0<\lambda<\lambda_0$, $C_0>1$, there exist
 $\,\xi=\xi(\lambda,v,C_0)>0,$ $\epsilon_0=\epsilon_0(v,\lambda)>0$,
 $\epsilon_1=\epsilon_1(v,\lambda,C_0)>0$,
such that if $\alpha \in \EDC(\nu,\xi)$, then the family
$\hat H_{\lambda v,\alpha,\theta}$ is almost localized
with parameters $C_0,\epsilon_0,\epsilon_1,C_1$, where
$C_1=C_1(v,\lambda,C_0,\nu)>0$.  For $v(x)=\pm 2 \cos 2 \pi x$,
$\lambda_0=1$.
\et

\subsection{Proof of Theorem \ref {SLE}}

To simplify the notation, fix some $v:\R/\Z \to \R$ analytic, and
set $\hH=\hH_{\alpha,\theta} =
\frac {1} {\lambda}\hat H_{\lambda v,\alpha,\theta}$.
Through this section, $C$ denotes an absolute large
constant, while $C_\sigma$, for instance, denote a
large constant that only depends on $\sigma$.  We warn that this is not the
same convention used in the remaining of this paper.

%A formal solution $\hat{u}_x $ of the equation
% $\hat H
%%_{\lambda v,\al,\theta}\hat{u} = E\hat{u}$ will
%be called {\it a generalized eigenfunction} if 
%$\la{ge} |\hat{u}_x| \le 1+|x|$
% The corresponding $E$ is called a {\it generalized eigenvalue.}

Fix an interval $I \subset \bbz.$ Let $\Gamma$ be the coupling operator
between $I$ and $\bbz\backslash I:$
$$ \Gamma(i,j)= \left\{\;
\begin{array}{cc}
& \hat v_{i-j},\; \chi_I(i)+\chi_I(j)=1,\\
& 0,\;\;\;\mbox{otherwise}.
\end{array} \right.$$

Then we can write $\hat{u}=-(\hH-E-\Gamma)^{-1}\Gamma\hat{u}$ from which
for any $x\in I$ we obtain 
 \be \la{poi}
 \hat{u}(x)=-\displaystyle\sum_{y\in
I,k\notin I} G_I (x,y)\hat v_{y-k}\hat{u}(k),
\ee
where $G_I$ is the Green's
function of $\hH$ restricted to the interval $I:$
$G_I=G_I(E)=(R_I(\hH-E)R^*_I)^{-1}$ (where $R_I:l^2(\Z) \to l^2(I)$ and
$R_I^*:l^2(I) \to l^2(\Z)$ denote the restriction and the inclusion
operators).
Set $a_k=\displaystyle\sum_{j: |j| \geq |k|,jk\ge 0} |j\hat v_j|.$
\comm{
By analyticity of $v$ we have
\be \la{anal}a_k\le Ce^{-\mu|k|}\ee with $C<\infty,\mu>0$
depending only on $v.$
}

\vskip .15in 
%\noindent{\bf Definition.}
Fix $E \in \bbr, \;m\in\bbr.$ A point $y \in {\bbz}$
 will be called 
$(m,k)$-regular if there exists an interval $I=[x_1+1,x_2-1],\;x_2=x_1+k+1,$
 such that $x_1<y<x_2$ and
$$\displaystyle\sum_{x\in I,i=1,2}|G_I(y,x)a_{x-x_i}|<e^{- mk}$$
%,\;\mbox{and}\;\;{\rm dist} (y,x_i) \ge{\frac {1} {5}} k;\; i=1,2. $$
 Otherwise, $y$ will be called $(m,k)$-singular.

%By (\ref{poi}) it suffices to
The strategy is to show that the existence of a generalized eigenfunction
as in Definition \ref {sle} implies that $k$ is $(m,\ell(k))$ regular
for an appropriate $m>0$ and for $\ell(k)$ comparable with $k,$ with $k$ in the desired ``between the resonances'' region.

%Similarly, if $k$ is $(m,k)$-regular then  $|$\hat u_k|$ is
%bounded by $Ce^{-mk}.$ hat we need is to establish 
%
% so it suffices to prove that for an appropriate
%$m$ and sufficiently large $k,$ $(m,k)$-singular points cannot be at a
%distance $k$ from, say, zero.

Define $I\subset \bbz$ by
$I=[0,N-1],$
%\,x\in\bbz,
$N\in\bbn.$ 
%Fix $E\in \bbr.$ 
We will omit the $\alpha, E$ dependence of
various quantities in what follows, and all constants
will be uniform for all $E$ in the spectrum.
We will also often assume $v\not\equiv 0,$ as
otherwise our statements become trivial.
 Let $P_N(\theta)=\det R_I(\hH_{\alpha,\theta}-E)R^*_I$. 
$P_N(\theta)$ is an even function of $\theta + \frac {N-1}{2}\alpha$ and
 can be written as a polynomial of degree $N$ in $\cos 2\pi(\theta+\frac {N-1}
{2}\alpha):$ $ P_N(\theta)
=\sum_{j=0}^N
c_j\cos ^j 2\pi(\theta + \frac {N-1}{2}\alpha)
\stackrel{\mathrm{def}}{=}Q_N(\cos 2\pi(\theta + \frac {N-1}{2}\alpha)).$

\bl \la{herman}

$\int_0^{1}{\frac {1} {N}}\ln |P_N(\theta)|d\theta \ge -\ln \lb.$

\el

\begin{pf}

It is proved by a standard Herman's subharmonicity argument \ci{H}.
\end{pf}

Let $A_{k,r}=\{\theta \in \R,\, |Q_k(\cos 2\pi\theta)|\le e^{(k+1)r}\}$.
The next lemma shows that every singular point ``produces" a long piece of the
trajectory of the rotation consisting of points 
belonging to an appropriate $A_{k,r}$.  It is rather immediate in the
almost Mathieu case, and we will adapt the argument of \ci{BJ1} for the general
case.

\bl \la{cramer}

There exists $\lambda_0=\lambda_0(v)>0$ such that for
$0<\lb<\lb_0,$ $\frac {1} {40} \leq \delta<\frac 12,$\footnote {As will be clear from the proof, the
result still holds for any $0<\delta<\frac {1} {2}$, but the constants,
including $\lambda_0$, become dependent on $\delta$.}
$\eps=\eps(\lb,v)>0$, $c=c(\lb,v)>0$, $K \subset \R \setminus \Q$
compact, $\alpha \in K,$ if
%$|P_N(\theta+x\alpha)|>(\lb^{-1}-\eps)^N,$  then all points $y\in x+I$ with 
%$dist (y, \partial (x+I))\ge \delta N$ are $(m,N)-$regular for sufficiently large $N.$.
$y\;\in \; {\bbz}$ is $(c,N)$-singular,
$N>N(\lambda,v,K)$ and $x \in \bbz$ is such that
$y-(1-\delta) N\le x \le y-\delta N,$
we have that $\theta+(x+ \frac{N-1}{2})\al$  belongs to $A_{N,-\ln
  \lb-\eps}.$ For $v(x)=\pm 2 \cos 2\pi x,\;\lb_0=1.$

\el

Assume, without loss of generality, that
$C_0|n_{j_k}|<k<\frac {1} {C_0} |n_{j_k+1}|$,
and $k$ is large (depending only on
$C_0$).  We will define scales $n \geq 0$ and $s \geq 1$ associated with
$k$ so that
\be \la{scales}
2s q_n \leq \zeta k<\min\{2(s+1)q_n,2q_{n+1}\},
\ee
where $\zeta=\frac {1} {32}$
if $2 |n_{j_k}|<k<\frac {1} {2} |n_{j_k+1}|$ and
$\zeta=\frac {C_0-1} {16 C_0}$ otherwise. Note that $s,n$ depend on $\epsilon_0.$

\bl \la{reg}

Assume there exists $\hat u$ as in Definition \ref {sle}.
For $\lb,c$ as in Lemma \ref{cramer}, there exist
$\xi=\xi(v,\lambda,C_0)>0$,
$\eps_0=\epsilon_0(v,\lambda)>0$, such that if $\alpha \in
\EDC(\nu,\xi)$,
$k>k(v,\lambda,c,C_0,\nu,\xi)$ and $s$, $n$ are as above
then $k$ is $(c,6sq_n-1)$-regular.

\el

Since  $sq_n>C_{C_0}^{-1} k$ the theorem now follows immediately
from the definition of regularity and (\ref{poi}). It therefore suffices to
prove Lemmas \ref{cramer}, \ref{reg}.

\subsubsection{Proof of Lemma \ref{cramer}}

Without loss, set $x=0.$ Set
 $x_1=-1,\;x_2=N.$ Assume $\theta+\frac {N-1} {2} \alpha
\notin A_{N,-\ln\lb-\eps}$, that is $P_N(\theta)>\lambda^N e^{-\epsilon N}$.
We need to show
that for $y\in [x_1,x_2] $ with $\mbox{dist} (y,\partial [x_1,x_2])
\ge \delta N$ we have
\be \label {star}
(*)=\sum_{z\in
I,\,i=1,2}|G_{I}(y,z)||a_{z-x_i}|<e^{-cN}.
\ee
By Cramer's rule $G_{I}(y,z)=\frac {\mu_{y,z}} {P_N(\theta)}$ where
$\mu_{y,z}$ is the corresponding minor.  The following lemma
reduces the study of Green's function to the study of determinants of
restrictions of $\hH-E$.

\bl [Lemma 10, \ci{BJ1}]

\be   \la{path}\mu_{y,z}=\displaystyle\sum_{\gamma}\alpha_{\gamma}\det
R_{I\backslash \gamma}(\hH-E)R^*_{I\backslash\gamma}
\displaystyle\prod_{i=1}^{|\gamma|}|\hat v_{\gamma_{i+1}-\gamma_i}|\ee
where the sum is taken over all ordered subsets
$\gamma=(\gamma_1,\ldots,\gamma_n)$ of $I$ with $\gamma_1=y$ and
$\gamma_n=z,\;|\gamma|=n-1,$
and  $\alpha_{\gamma}\in\{-1,1\}.$

\el
%{\bf Proof.} Expand $\mu_{y,z}$ through the $y-$column:
%$\mu_{y,z}=\displaystyle\sum_{\gamma_2\not= y}\delta_{\gamma_2}
% a_{\gamma_2-y}h_{z,y,\gamma_2}$
%where $\delta_{\gamma_2}\in\{-1,1\},$ and $h_{z,y,\gamma_2}$ is the determinant of the
%original matrix with $z$ and $y$ raws and $y$ and $\gamma_2$ columns removed.
%Continue this process with each term until $\gamma_n=z.$ The
%contribution of each term is equal, up to a sign, to  $\det
%R_{I\backslash \gamma}(H-E)R_{I\backslash\gamma}
%\displaystyle\prod_{i=1}^{|\gamma|}|\hat v_
%{\gamma_{i+1}-\gamma_i}|.$\qed

%Note that already an obvious bound of the form 

%\be \la{dur}
%|\det R_{I\backslash \Lambda}(H-E)R_{I\backslash\Lambda}|\le
% (2|\lb|+|E|+\| v\|_\infty)^{|I\backslash \Lambda|}
%\ee
%would provide a result for sufficiently large $\mu,$ however a more delicate estimate is needed for small $\mu.$ Assume $\omega$ is Diophantine.

The following is a simplified, yet more general,
version of Lemma 11 of \cite {BJ1}.  It gives an upper bound on
$\mu_{y,z}$.

\bl \la{upper}

For any $\Lambda \subset I,$
for sufficiently large $N>N(v,\lambda,K)$,
\be  \la{1}
|\det R_{I\backslash \Lambda}(\hH-E)R^*_{I\backslash\Lambda}|\le
\lambda^{-N} e^{C \|v\|_0^{1/2} \lambda^{1/2} N}
\left (\|v\|_0+C^{-1} \lambda^{-1}
\frac {\# \Lambda^2} {N^2} \right )^{-\#\Lambda}.
\ee
\comm{
and if $\#\Lambda>\lb^{1/10}N$ then 
\be \la{2}
|\det R_{I\backslash \Lambda}(\hH-E)R_{I\backslash\Lambda}|\le \lb^{-(N-1/5|\Lambda|)}  
\ee
}

\el
\begin{rem}
Unlike the corresponding upper bound in \cite{BJ1} this Lemma does not require $\alpha \in \DC.$ 
\end{rem}  
\begin{pf}

By Hadamard's bound, we have
\be
|\det R_{I\backslash \Lambda}(\hH-E)R^*_{I\backslash\Lambda}| \leq
\prod_{j \in I \setminus \Lambda}
((2\lambda^{-1} \cos 2 \pi (\theta+j\alpha)-E)^2+\|v\|_{L^2}^2)^{1/2}.
\ee
Thus
\be
\ln |\det R_{I\backslash \Lambda}(\hH-E)R^*_{I\backslash\Lambda}| \leq
\sum_{j \in I \setminus \Lambda}
\frac {1} {2} \ln
((2\lambda^{-1} \cos 2 \pi (\theta+j\alpha)-E)^2+\|v\|_0^2)=\sum_{j \in I \setminus
\Lambda} u(x+j\alpha).
\ee
Set $A_t= \{x:\, |2 \cos (2 \pi x)-\lambda E|<t\}$, $0 \leq t \leq
2+|\lambda E|$, and let $t_0$ be such that the Lebesque
measure of $A_{t_0}$ is equal to
$\frac {\# \Lambda}N.$ Then by unique ergodicity of
irrational rotations and continuity of $u$, we have
\be
\sum_{j \in I \setminus
\Lambda} u(x+j\alpha) \leq \sum_{j \in I } u(x+j\alpha)-\sum_{\ntop {j \in
I,} {x+j\alpha \in A_{t_0}}}
u(x+j\alpha)+N o(1)\leq N(\int u(x)(1-\chi_{A_{t_0}}) dx +o(1)),
\ee

%\be
%\sum_{j \in I \setminus
%\Lambda} u(x+j\alpha) \leq N (\int u(x) d\omega(x)+o(1)),
%\ee
%for some measure $\omega$, absolutely continuous with density bounded by
%$1$, of mass $1-\frac {|\Lambda|} {N}$.
Notice that a direct computation gives
\be
\int u(x) dx=\ln \lambda^{-1}+|\ln |z+\sqrt{z^2-1}||,
\ee
where $z=\frac {E+i \|v\|_0} {2 \lambda^{-1}}$ (see \cite {BJ1} for
details).  From
this and $|E| \leq 2\lambda^{-1}+\|v\|_0$ in the spectrum,
we get the estimate
$\int u(x) dx \leq \ln \lambda^{-1}+C \lambda^{1/2} \|v\|_0^{1/2}$.
To complete the proof, it suffices to show that, letting
$s=\frac {\# \Lambda} {N}$,
\be \label {ln 1+C}
\int u(x) \chi_{A_{t_0}}(x) dx \geq
s \ln (\|v\|_0+C^{-1} \lambda^{-1} s^2)-
C \|v\|^{1/2}_0 \lambda^{1/2}.
\ee
%when $\phi(x)= 1-\chi_B(x)$, with $\int \phi(x) dx=t$.
%This is clear when $s^2<\lambda \|v\|_0$, so we assume $s^2 \geq \lambda
%\|v\|_0$ from now on.

Let us split $A_t$ into four segments $I_j$ of length $s_j \leq 1$
such that $x \mapsto (2 \cos (2 \pi x)-\lambda E)^2$ is monotonic
in $I_j$.  Then, since $1-\cos 2\pi x \geq 16 x^2$ for $x\leq 1/4,$
\begin{align}
\int_{I_j} u(x) dx &\geq \int_{0}^{s_j} \frac {1} {2}
\ln \left (\|v\|_0^2+4 \lambda^{-2} (1-\cos (2 \pi x))^2 \right ) dx \geq
\int_{0}^{s_j} \frac {1} {2} \ln \left (\|v\|_0^2+256 \lambda^{-2}
x^4 \right ) dx\\
\nonumber
&\geq \frac {1} {2} s_j
\ln \left (\|v\|^2_0+256 \lambda^{-2} s_j^4 \right )-2 s_j.
\end{align}
Since $x \mapsto x \ln \left (\|v\|^2_0+
256 \lambda^{-2} x^4 \right )$ is concave on $\R^+$, we get, summing over
$j$,
\be
\int_{A_{t_0}} u(x) dx \geq \frac {1} {2} s \ln \left (\|v\|_0^2+
\lambda^{-2} s^4 \right )-2s.
\ee
Considering separately the cases $\lambda^{-1}s^2>c\|v\|_0$ and
$\lambda^{-1}s^2<c\|v\|_0$ we obtain that this implies (\ref {ln 1+C}).
\end{pf}

Our assumption $|P_N(\theta)|>\lb^{-N}e^{-\eps N},$ implies, with the
notation of (\ref{star}) and
using (\ref{path}):
\begin{align} \label {sums}
(*) &\le (\lb e^{\eps})^{N} \sum_{n=1}^{N-1}
\sum_{\ntop{i=1,2,} {\gamma:\,|\gamma|=n}}|\det
R_{I\backslash \gamma}(\hH-E)R^*_{I\backslash\gamma}|
|a_{x_i-\gamma_{|\gamma|+1}}|
\prod_{i=1}^{n}|\hat v_{\gamma_{i+1}-\gamma_i}|\\
\nonumber
&\leq e^{(\epsilon+C \|v\|_0^{1/2} \lambda^{1/2}) N} \sum_{n=1}^{N-1}
\sum_{\ntop{i=1,2,} {\gamma:\,|\gamma|=n}} C_\sigma C_{v,\sigma}^{n+1}
\left (\|v\|_0+C^{-1} \lambda^{-1} \frac {(n+1)^2} {N^2} \right )^{-(n+1)}
e^{-\sigma b(\gamma,i)},
\end{align}
where $0<\sigma<\sigma(v) \leq 1$ is such that
\be\label{cond1}|\hat v_k| \leq C_{v,\sigma} e^{-|k| \sigma},\ee and
$b(\gamma,i)=|\gamma_{|\gamma|+1}-x_i|+\sum_{i=1}^{|\gamma|}
|\gamma_{i+1}-\gamma_i|$.
Let $G_{b,n}=\{\gamma,\, |\gamma|=n \text { and } b(\gamma,i)=b\}$.  Then
\begin{align}
(*)&\le
e^{(\epsilon+C \|v\|_0^{1/2} \lambda^{1/2}) N} \sum_{n=1}^{N-1} \sum_b
C_\sigma C^{n+1}_{v,\sigma}
\left (\|v\|_0+C^{-1} \lambda^{-1} \frac {(n+1)^2} {N^2} \right )^{-(n+1)}
e^{-\sigma b} \# G_{b,n}\\
\nonumber
&\leq
e^{(\epsilon+C \|v\|_0^{1/2} \lambda^{1/2}) N} \sum_{n=1}^{N-1}
C_\sigma (2 C_{v,\sigma})^{n+1}
\left (C^{-1} \lambda^{-1} \frac {(n+1)^2} {N^2} \right )^{-(n+1)}
\sum_{b,\, G_{b,n} \neq \emptyset}
e^{-\sigma b} \binom{b}{n}.
\end{align}
If $G_{b,n} \neq \emptyset$ then $\delta N\leq \max \{\dist(y,\partial I),n+1\}
\leq b \leq (n+1) N\leq N^2$.
By Stirling formula, setting $b=r N$, $n+1=s b$, we have $\binom {b} {n} \leq
C r N e^{\phi(s) r N}$ where $\phi(s)=-s \ln
s-(1-s)\ln(1-s)$.  Thus we can estimate
\be
(*)\le
e^{(\epsilon+C \|v\|_0^{1/2} \lambda^{1/2}) N} C C_\sigma N^5
\sup_{\ntop {\delta \leq r \leq n+1} {0<s \leq 1}}
\left (\frac {C^{-1} \lambda^{-1}} {2 C_{v,\sigma}}
r^2 s^2 \right )^{-r s N}
e^{-\sigma r N} e^{\phi(s) r N}
\ee
The desired exponential decay of $(*)$ follows if
\be
(**)=\sup_{0<s \leq 1}
\epsilon+C\|v\|_0^{1/2} \lambda^{1/2}+\left
(\ln 2 C+\ln C_{v,\sigma}+\ln \lambda
-2 \ln \delta s-\frac {\sigma} {s}+\frac {\phi(s)} {s} \right ) \delta s<0.
\ee
This condition is satisfied (for appropriate $\epsilon$)
if \be \label{cond2} c_0=\delta^{-2} \lambda C_{v,\sigma} \sigma^{-3}\ee is small.  Indeed,
using that $\|v\|_0 \leq \frac {C C_{v,\sigma}} {\sigma}$, and that
$\phi(s)/s \leq 1 -\ln s$,
\be \label{518}
(**) \leq \epsilon+
\left (C \delta c_0^{1/2}-\frac {\delta} {2} \right ) \sigma+
\left (\ln C+\ln c_0+3 \ln \frac {\sigma} {s}-
\frac {\sigma} {2 s} \right ) \delta s.
\ee

\comm{
Using (\ref{1}), we bound $\textstyle\sum_1$ by
$$
\textstyle\sum_1\le\lb^{-N}e^{C_v\lb^{\frac {1}
{2}}N}\displaystyle\sum_{n=1}^{\lb^{\frac {1} {10}}N}C^n \displaystyle\sum_{i=1,2;\,b}
e^{-\mu b}
%\times$$
%$$\times 
\#\{\gamma:\,
|\gamma|=n,\,\gamma_1=y,\,|\gamma_{n+1}-x_i|+
\displaystyle\sum_{i=1}^n|\gamma_{i+1}-\gamma_i|=b\}$$
$$\le\lb^{-N}e^{C_v\lb^{\frac {1} {2}}N} 
\displaystyle\sum_{b\ge \dist(y,\partial I)}e^{-\mu b}
\displaystyle\sum_{n=1}^{\lb^{\frac {1} {10}}N}C^n \binom{b}{n}$$
where $C,\mu$ are from \ref{anal}.
Assume $C\ge 1$ (otherwise bound it by $1.$) For $y$ with $\mbox{dist} (y,\partial [x_1,x_2])\ge \delta N$  
we have, for small $\lb$
$$
\textstyle\sum_1\le \lb^{-N}e^{C_v\lb^{\frac {1} {2}}N}\lb^{\frac {1} {10}}N
C^{\lb^{\frac {1} {10}}N}\displaystyle\sum_{b\ge \dist(y,\partial
I)}e^{-\mu b}\binom{b}{\lb^{\frac {1} {10}}N}.
$$
By Stierling formula, for large $N$,
$$
\binom{b}{\lb^{\frac {1} {10}}N}
\le \binom {b}{\delta^{-1}\lb^{\frac {1} {10}}b}\le
C_\lambda e^{\phi(\delta^{-1}\lb^{\frac {1} {10}})b}
$$
where $\phi(x)=-x\log x -(1-x)\log x.$

Hence,
\be \la{sum1}
\textstyle\sum_1 \le \lb^{-N} C_{\lambda,v}
e^{(C_v\lb^{\frac {1} {2}}+
\lb^{\frac {1} {10}}\log C_v +\delta\phi(\delta^{-1}
\lb^{\frac {1} {10}})-{\mu \delta})N}.
\ee
Similarly, using (\ref{2}), we have 
\beq \la{sum2}
\textstyle\sum_2 &\le& \lb^{-N}e^{C_v\lb^{\frac {1} {2}}N}\displaystyle\sum_n C^n
\lb^{\frac {1} {5}n}%\nonumber
\\&\times&\displaystyle\sum_{b\ge \dist(y,\partial I)}e^{-\mu
b}\binom{b}{n} \le  \lb^{-N}e^{C_v\lb^{\frac {1} {2}}N}\displaystyle\sum_
{b\ge \dist(y,\partial I)}e^{-\mu b}
(1+C_v\lb^{-\frac {1} {5}})^b\nonumber
%\nonumber
\eeq

Substituting (\ref{sum1}),(\ref{sum2}) into (\ref{sums}) we obtain the
desired exponential decay for sufficiently small $\epsilon$ and large $N$ provided
\be \la{est3}
\eps+C_v\lb^{\frac {1} {2}}+
\lb^{\frac {1} {10}} C_v+\delta\phi(\delta^{-1}
\lb^{\frac {1} {10}})-{\mu \delta}<0\ee
and
\be\la{est4}
\eps+C_v\lb^{\frac {1} {2}}+\delta\log (1+C_v \lb^{\frac {1} {5}})-{\mu\delta}<0\ee
}

In case $v(x)=2 \cos 2 \pi x$ the result in the form we need, for any
$0<\lb<1$ is stated in Lemma 9.2  in \ci{AJ} (uniformity on $E$
was not claimed but is automatic from the argument).
The simplification in this case is mainly due to the fact
that, as can be seen by a direct computation (but also follows immediately
from Lemma~\ref{path}),
$|\mu_{y,z}|=|P_{y-1}(\theta)P_{N-z}(\theta+z\alpha)|$
and we have a uniform upper bound of the form
$\ln |P_k (\theta)|<C_{\epsilon,K}+k(\ln \lambda^{-1}+\epsilon)$.
It is obtained by identifying
$P_k(\theta)$ with the upper-left coefficient of the
$k$-th iterate of the almost Mathieu cocycle
$(\alpha,S_{\lambda^{-1}v,E})$, and applying
the formula $L(\alpha,S_{\lambda^{-1}v,E})=-\ln \lambda$ for $0<\lambda<1$
and all $E$ in the spectrum \cite {BJ2}.
\qed

%In order to obtain $\lb_0$ close to
%$1$ for functions $f$ close to the cosine one needs to consider
%separately terms with $\gamma$ such that $I\backslash \gamma$ contains
%large (of size $\approx N$) intervals, replace estimate from Lemma \ref{upper}
%by a uniform upper bound with $\int\ln
%P_k(\th)+\epsilon$ (see \ci{j} or \ci{bg}) for such $\gamma,$ and
%choose
%parameters so that $|a_k|(4|\lb^{-1}|+2\| f\|_\infty)^k <1$ for  $k\ge 2.$ 
%

%We will show that $k$ is $(L,3sq_n-1)-$regular. From this the statement follows immediately.

\subsubsection{Proof of Lemma \ref{reg}}
We first recall some concepts introduced in \ci{AJ}, Sec. 9.
\begin{definition}

We will say that the set $\{\theta_1,\ldots,\theta_{k+1}\}$ is
{\it $\eps$-uniform} if 
%for any $z \in [-1,1]$ 

\be \la{prod}
\max_{z\in[-1,1]}\max_{j=1,\ldots,k+1}\prod_{^{\ell=1}_{\ell \not= j}}^{k+1}\frac{
|z-\cos 2\pi\theta_\ell)|}
{|\cos 2\pi\theta_j-\cos 2\pi\theta_\ell)|}
%\stackrel{\mathrm
%{def}}{=}\frac{|I_1|}{|I_2|} 
< e^{{k\eps}}
\ee

\end{definition}

 $\eps$-uniformity (the smaller $\eps$ the better) involves uniformity along with certain cumulative
repulsion
of $\pm \theta_i$(mod $1$)'s.
%For example, as will follow from the estimates
%below, the set $\{\frac{i}{k+1}\}_{i=0,\ldots,k}$ is $\eps$-uniform
%for any $\eps >0.$
%This is not true, as $\cos 2\pi/(k+1)=\cos 2 \pi k/(k+1)$. 
%good point!

\bl \la{nongen}
Let $\hat \eps>0.$ If $\theta_1,\ldots,\theta_{k+1}\in A_{k,-\ln \lb-\eps}$ and
$k>k(\eps,\hat \eps)$ is sufficiently large, then
$\{\theta_1,\ldots,\theta_{k+1}\}$ is not $(\epsilon-\hat \eps)$-uniform.
\el

\begin{pf}

Write polynomial $Q_k(z)$
 in the Lagrange interpolation form using $\cos 2\pi\theta_1,\ldots,\cos
2\pi\theta_{k+1}:$
\be \la{lagr}
|Q_k(z)|=\left
|\sum_{j=1}^{k+1}Q_k(\cos 2\pi\theta_j)\frac {\prod_{\ell \not= j}
(z-\cos 2\pi\theta_\ell)}
{\prod_{\ell \not= j}(\cos2\pi\theta_j-\cos 2\pi\theta_\ell)} \right |
\ee
 Let $\theta_0$ be such
 that $|P_k(\theta_0)|\ge
 \lb^{-1}$.  The lemma now follows immediately from
 (\ref{lagr}) with $z=\cos 2\pi(\theta_0+\frac {k-1}{2}\al)$.
\end{pf}

%We will first prove the theorem for sufficiently large $j,$ so that $|n_{j-1}|>C
%(\kappa,\tau),$ and with $C_0=3,$ the theorem as stated will then follow by adjusting 
%$C_0$ depending on $\eps_0$ and the value of $C(\kappa,\tau).$

%Let $p=\max\{i<j: n_i\le 0\}.$ 
We now define $I_1,I_2\subset \bbz$ as follows
\begin{enumerate}
\item $I_1=[-2sq_n+1,0]$
% if $\ell_0> 0,$ 
and $I_2=[k-2sq_n+1,k+2sq_n],$ if $k< \frac {|n_{j_k+1}|}3$ and $n_{j_k} \geq 0$.
\item $I_1=[1,2sq_n]$
% if $\ell_0> 0,$ 
and $I_2=[k-2sq_n+1,k+2sq_n],$ if $k< \frac {|n_{j_k+1}|}3$ and $n_{j_k}<0$.
\item $I_1=[-2sq_n+1
%+n_p,n_p
,2sq_n]$ and $I_2=[k-2sq_n+1,k],$ if  $\frac {|n_{j_k+1}|}3\le k<
\frac {|n_{j_k+1}|}2$. 
\item  $I_1=[-2sq_n+1
%+n_p,n_p
,2sq_n]$ and $I_2=[k+1,k+2sq_n],$ if  $k \geq \frac {|n_{j_k+1|}}2$.
\end{enumerate} 
In either case, the total number of elements in $I_1\cup I_2$ is $6sq_n.$ Set $\theta_j=\theta+j\alpha,j\in I_1\cup I_2.$ 

\bl \la{uni}
 If $\xi<C^{-1} \epsilon_0$ and
$k>k(C_0,\nu,\xi)$ then set $\{\theta_j\}_{j\in
I_1\cup I_2}$ is $C \eps_0+C_{C_0} \xi$-uniform.
\el

\begin{pf}

We will first estimate the numerator in (\ref{prod}). We have, in each case,

\begin{align} \la{num}
\sum_{^{j\in I_1\cup I_2}_{j\not=i}}&\ln |\cos 2\pi a-\cos
2\pi \theta_j|\\
\nonumber
&=\sum_{^{j\in I_1\cup I_2}_{j\not=i}}\ln
|\sin 2\pi\frac{a+\theta_j}2|+\displaystyle\sum_{^{j\in I_1\cup
    I_2}_{j\not=i}}\ln |\sin 2\pi\frac{a-\theta_j}2|+(6sq_{n}-1)\ln
2\\
\nonumber
&=\Sigma_+ + \Sigma_- + (6sq_{n}-1)\ln 2.
\end{align}

Both $\Sigma_+$ and $\Sigma_-$ consist of $6s$ terms
of the form of (\ref{irrat-1}) plus $6s$ terms of the form
\be
\ln \min_{j=1,\ldots,q_{n-1}}|\sin( 2\pi(x+\frac{j\al}{2}))|,
\ee
minus $\ln |\sin\frac{a\pm\theta_i}2|$. 
%It is easily seen that
%\be \la{min}
%\min_{i=1,\ldots,q_{n-1}}|\sin(
%2\pi(x+\frac{k\al}{2})|<
%\frac{\pi}2(\Delta_{n-1}+\Delta_{n-2})<\frac{\pi}{q_{n-1}}
%\ee
Therefore, by (\ref{irrat-1})

\be\la{num1}\displaystyle\sum_{^{j\in I_1\cup I_2}_{j\not=i}}\ln |\cos 2\pi a-
\cos
2\pi\theta_j|\le -6sq_{n}\ln 2 + Cs\ln q_{n}.
%2s \ln 2\pi + o(s)
\ee

To estimate the denominator of (\ref{prod}) we represent it again in
the form (\ref{num}) with $a= \theta_i.$
 Then
\be
\Sigma_-=\sum_{^{j\in I_1\cup I_2}_{j\not=i}}\ln |\sin \pi
(i-j)\al|.
\ee
$\Sigma_-$ consists of $6s$ terms
of the form of (\ref{irrat-1}) plus $6s-1$ minimum terms (since, when splited 
into $6s$ sums over intervals of length $q_n$ each, one of the sums will be 
exactly of the form of one of (\ref{irrat-1})). For each $i,j\in
I_1\cup I_2,$ we have $|i-j|\le k+4sq_n<C_{C_0}sq_n$.
In particular $\ln |\sin \pi (i-j)\al|\ge -C_{C_0} s q_n \xi,$ for
large $q_n$ (depending on $\nu,\xi$).  Notice that
\be \label {alphaxj}
\max \{\ln |\sin x|,\ln |\sin x+\pi j\alpha|\}>-2 \xi q_n,\quad x \in
\R,0<|j|<q_{n+1},
\ee
provided $q_n$ is sufficiently large (depending on $\nu$ and $\xi$).  Using
that $s q_n<q_{n+1}$ we get
 \be\la{den1}
\Sigma_-\ge -6sq_{n}\ln 2 - C_{C_0} s q_n \xi.
\ee
Similarly, $\Sigma_+$ consists of $6s$ terms
of the form of (\ref{irrat-1}) plus $6s$ minimum terms, each of the form 
\be
\ln 
%\min_{j=1,\ldots,q_{n-1}}
|\sin( 2\pi(\theta+\frac{(i+j)\al}{2}))|
\ee
for some $|j|<C_{C_0} s q_n<C_{C_0} q_{n+1}$,
minus $\ln |\sin 2 \pi (\theta+i\alpha)|$ (which cancels a possible minimal
term with $i=j$).
Using (\ref {alphaxj}), we see that at most $6$ minimum
terms are smaller than $-2 \xi q_n$.   Let us estimate the smallest term
with $j \neq i$.

Consider first cases 1 and 2 of the definition of $I_1$, $I_2$.
Then, by the definition of
$I_i,\,i=1,2,$ and using (\ref{scales}) we have  $i+j \neq n_{j_k}$
and $|i+j|<|n_{j_k+1}|$.  Therefore the smallest term
is bounded below by $-C \epsilon_0 s q_n$
when $k>2 |n_{j_k}|$, and by $-C_{C_0} \xi s q_n$, when $k \leq 2 |n_{j_k}|$.
In cases 3 and 4 we have $i+j \neq n_{j_k+1}$, $|i+j|<2 |n_{j_k+1}|$,
and $sq_n \geq C_{C_0}^{-1} |n_{j_k+1}|$.  Therefore the smallest term
is bounded below by 
$-C \xi |n_{j_k+1}|>-C_{C_0} \xi sq_n$.

Putting all together, we see that
\be\la{den2}
\Sigma_+\ge -6sq_{n}\ln 2 - Cs\ln q_{n}-(C \eps_0+C_{C_0}\xi) sq_n.
\ee 

%that either $|i+j|\le 2sq_n$ or 
%$n_j< k-sq_n\le i+j \le k+2sq_n < n_{j+1}.$ 

Combining, (\ref{num1}),(\ref{den1}), and (\ref{den2}), we obtain
$$
\max_{j\in I_1\cup I_2}\prod_{^{\ell\in I_1\cup I_2}_{\ell \not= j}}\frac{
|z-\cos 2\pi\theta_\ell|}
{|\cos 2\pi\theta_j-\cos 2\pi\theta_\ell|}<
e^{(C \epsilon_0+C_{C_0} \xi) sq_n+Cs\ln q_{n}}
$$
as desired.
\end{pf}

We can now finish the proof of Lemma~\ref{reg}.
By Lemmas \ref{nongen} and \ref{uni} at least one
 of $\theta_j,\;j\in I_1\cup I_2,$ is not in
 $A_{6sq_n-1,-\ln\lb-(C_{C_0}\xi+C \eps_0)-\hat \epsilon}$ where 
$\hat \epsilon$
can be made arbitrarily small for large $n$.
From (\ref{poi}) we see that the existence of a generalized eigenfunction
$\hat u$ (as in Definition (\ref {sle})) implies that $y=0$
is $(m,k)$-singular for any $m>0$ and for $k$ sufficiently large depending
on $m$.
Thus by Lemma \ref{cramer} with sufficiently small $\eps_0,$ $\xi$
(depending also on $C_0$), $y=0$ and
$\delta<\frac 16,$
%and singularity
%of $0$,\footnote {To get what we need here one can take in
%Lemma \ref {prop}, besides $y=0$, also
%$\epsilon=\frac {99} {100} L$ and $\delta=\frac {99} {400}$.}
we have that for all $j\in I_1,\;$  $\theta_j\in
A_{6sq_n-1,-\ln\lb-(C_{C_0}\xi+C \eps_0)-\hat \epsilon}$ for
sufficiently large $n.$  Let $j_0\in I_2$ be such
 that $\theta_{j_0}\notin A_{6sq_n-1,-\ln\lb-(C_{C_0} \xi+C \eps_0)-\hat
\epsilon}.$
Since $|j_0-k|\le 2sq_n$ we can apply again Lemma~\ref{cramer}
with $\delta<\frac 16$ to obtain that 
$k$ is $(c,6sq_n-1)$-regular.\qed

\section{Preliminary dynamical estimates} \label {growth}

\subsection{Lagrange interpolation argument}

Given Fourier coefficients $\hat w=(\hat w_k)_{k \in \Z}$ and an interval $I
\subset \Z$, we let $w^I=\sum_{k \in I} \hat w_k e^{2 \pi i k x}$.  The
length of the interval $I=[a,b]$ is $|I|=b-a$.

We will say that a trigonometrical polynomial $p:\R/\Z \to \C$
has essential degree at most
$k$ if its Fourier coefficients outside
an interval $I$ of length $k$ are vanishing.

\begin{thm} \label {trig pol}

Let $1 \leq r \leq [q_{n+1}/q_n]$.
If $p$ has essential degree at most $k=r q_n-1$ and $x_0 \in \R/\Z$ then for
some absolute constant $C$
\be \label {d1}
\|p\|_0 \leq C q_{n+1}^{C r} \sup_{0 \leq j \leq k}
|p(x_0+j\alpha)|.
\ee

\end{thm}

\begin{pf}

We may assume that $p(x)=P(e^{2 \pi i x})$ where $P$ is a polynomial of
degree $k$.  Then by Lagrange interpolation,
\be
p(x)=\sum_{j=0}^k p(x_0+j\alpha) \prod_{\ntop {0 \leq l \leq k,} {l \neq j}}
\frac {e^{2 \pi i x}-e^{2 \pi i (x_0+l \alpha)}} {e^{2 \pi i
(x_0+j\alpha)}-e^{2 \pi i (x_0+l \alpha)}}.
\ee
Thus
\be
\ln \|p\|_0 \leq \ln r q_n+\ln \sup_{0 \leq j \leq k}
|p(x+j\alpha)|+\sup_{\ntop {0 \leq j \leq k,}{x \in \R}}
\sum_{\ntop {0 \leq l \leq k,}{l \neq j}}
\ln \frac {|1-e^{2 \pi i (x+l \alpha)}|} {|1-e^{2 \pi i
(-j\alpha+l\alpha)}|}.
\ee

In order to prove (\ref {d1}) it is enough to
show that for $1 \leq s \leq r$, $0 \leq j \leq r
q_n-1$ and $x \in \R$ we have
\be \label {first est}
\sum_{\ntop {(s-1) q_n \leq l \leq s q_n-1,} {l \neq j}}
\ln |1-e^{2 \pi i (x+l \alpha)}| \leq C \ln q_{n},
\ee
\be \label {second est}
\sum_{\ntop {(s-1) q_n \leq l \leq s q_n-1,} {l \neq j}}
\ln |1-e^{2 \pi i (-j\alpha+l\alpha)}| \geq -C \ln q_{n+1}.
\ee

Consider first the case $(s-1)q_n \leq j \leq s q_n-1$.  Then
\be
\sum_{\ntop {(s-1) q_n \leq l \leq s q_n-1,} {l \neq j}}
\ln |1-e^{2 \pi i (x+l \alpha)}| \leq
\sum_{\ntop {(s-1) q_n \leq l \leq s q_n-1,} {l \neq l_0}}
\ln |1-e^{2 \pi i (x+l \alpha)}|,
\ee
where $|1-e^{2 \pi i (x+l_0 \alpha)}|$ is minimal.
Using that $|1-e^{2 \pi i y}|=2 |\sin \pi y|$ and Lemma \ref {9.7}, we get
(\ref {first est}).  The same argument gives (\ref {second est}) even
more directly, with $-C \ln q_n$ in the right hand side, since the sum to be estimated is already of the form
considered in Lemma \ref {9.7}.

Consider now the case $s-1 \neq [j/q_n]$.  Writing
\be
\sum_{\ntop {(s-1) q_n \leq l \leq s q_n-1,} {l \neq j}}
\ln |1-e^{2 \pi i (x+l \alpha)}|=
\sum_{(s-1) q_n \leq l \leq s q_n-1}
\ln |1-e^{2 \pi i (x+l \alpha)}|+\ln |1-e^{2 \pi i (x+l_0\alpha)}|,
\ee
where $(s-1) q_n \leq l_0 \leq s q_n-1$ is such that
$|1-e^{2 \pi i (x+l_0\alpha)}|$ is minimal, we see that
(\ref {first est})
follows from Lemma \ref {9.7}.  To obtain (\ref {second est}) from
Lemma \ref {9.7}, we must also show that
\be
\inf_{(s-1) q_n \leq l \leq s
q_n-1} \ln |1-e^{2 \pi i (l-j)\alpha}| \geq -C \ln q_{n+1}.
\ee
But this follows
from (\ref {b1}) and (\ref {b2}).
This concludes the proof of (\ref {d1}).
\end{pf}

The Diophantine condition $\alpha \in \DC(\kappa,\tau)$ implies
\be \label {b3}
q_{n+1} \leq \kappa^{-1} q_n^{\tau-1}.
\ee
In this case, (\ref {d1}) implies, with $c=\frac {1} {\tau-1}>0$ and
$C=C(\kappa,\tau)$,
\be \label {d2}
\|p\|_0 \leq C e^{C (1+k)^{1-c} \ln (1+k)} \sup_{0 \leq j \leq k}
|p(x+j\alpha)|,
\ee
since $1+k=r q_n \leq q_{n+1} \leq \kappa^{-1} q_n^{\tau-1}$.

\subsection{Polynomial growth}

From now on 

$\alpha \in \DC(\kappa,\tau)$ and the family $\{\hat H_{\lambda
v,\alpha,\theta}\}_{\theta \in \R}$ is almost localized with parameters $(\epsilon_0,C_0,C_1,\epsilon_1)$.
Let $E$ be in the spectrum and let $A=S_{\lambda v,E}$.  Let
$\theta=\theta(E)$ be given by Theorem \ref {boundedeigenfunctions},
and let $\{n_j\}$ be the set of resonances of $\theta(E)$.

All constants may depend
on $\kappa$, $\tau$, $C_0$, $\epsilon_0$,
$\epsilon_1$, $C_1$, and on bounds on the analytic extension of $\lambda
v$.  We will use $C$ to denote large constants and $c$ to denote
small constants.  Further dependence on other parameters,
will be explicitly indicated, for instance, we will use $C_\delta$ for a
large constant that depends on all parameters
above and also on an arbitrary parameter $\delta>0$.

The specific values of generic constants such as $C$, $c$, $C_\delta$, may
change through the argument, even when they appear in the same formula.  In a few places we will use
non-generic constants, denoted
$C^{(1)}$, $C^{(2)}$,..., to simplify later referencing.

\begin{thm}

We have $L(\alpha,A)=0$.

\end{thm}

\begin{pf}

Fixing a non-resonant $\theta'$, Theorem \ref {almost localized} (see Remark
\ref {localization non-resonant}) implies localization for
$\hat H_{\lambda v,\alpha,\theta'}$.
By Theorem \ref {localreduc}, $L(\alpha,S_{\lambda v,E'})=0$ for a dense
set of $E'$ in the spectrum.  By continuity of the
Lyapunov exponent, \cite {BJ2}, this implies that
$L(\alpha,S_{\lambda v,E})=0$ for all $E$ in the spectrum.\footnote{An alternative
argument for zero Lyapunov exponent to be dense is the
following.  By duality \cite{gjls}, if $\hat H$ has pure point spectrum for almost
every $\theta$ then $H$ has absolutely continuous spectrum for almost every
$\theta$ (see discussion in \cite {BJ1} for the models considered here).
By the Ishii-Pastur Theorem (e.g.,\cite {Cy}), the Lyapunov
exponent is zero densely in the spectrum.}
\end{pf}

In our context (see Remark \ref {uniform growth}), this means that
\be \label {context}
\sup_{x \in \R/\Z} \|A_k(x)\|
\leq C_\delta e^{\delta k}, \quad \delta>0.
\ee

Our goal in this section is to improve subexponential growth to polynomial
growth.  This is of course still not optimal
(see Corollary \ref {linear}), but it is a good starting point.

\begin{thm} \label {growth polynomial}

We have $\|A_s\|_0 \leq C (1+s)^C$, $s \geq 0$.

\end{thm}

We note that the proof of this estimate involves some themes which will
appear again later.

Choose $4 C_0 (|n_j|+1)<m<C_0^{-1} |n_{j+1}|$ of the form
$m=r q_k-1<q_{k+1}$, let $I=[-[m/2],m-[m/2]]$
and define $u(x)=u^I(x)$.  Let
$U(x)=\left (\bm e^{2\pi i \theta} u(x)\\u(x-\alpha) \em
\right )$.  Then
\be \la{ais}
A(x) \cdot U(x)-e^{2 \pi i \theta} U(x+\alpha)=e^{4 \pi i \theta}
\left (\bm h(x)\\0 \em \right ),
\ee
where
\be 
\hat h_k=\chi_I(k) 2 \cos 2 \pi (\theta+k \alpha) \hat u_k+\sum \chi_I(k-j)
\hat v_j \hat u_{k-j},
\ee
where $\chi_I$ is the characteristic function of $I$.
Since $\hat H \hat u=E \hat u$, we also have
\be
-\hat h_k=\chi_{\Z \setminus I}(k) 2 \cos 2 \pi (\theta+k \alpha)
\hat u_k+\sum \chi_{\Z \setminus I}(k-j) \hat v_j \hat u_{k-j}.
\ee
The estimates $|\hat u_k|<C_1 e^{-\epsilon_1 |k|}$ for
$m/4<|k|<m$, $|\hat u_k| \leq 1$ for all $k$ and
$|\hat v_k| \leq C e^{-c |k|}$ for all $k$ then
imply that $|\hat h_k| \leq C e^{-c m} e^{-c |k|}$, that is \be\la{hc}\|h\|_c
\leq C e^{-c m}.\ee

\begin{thm} \label {lower bound ux}

We have $\inf_{x \in \R/\Z} \|U(x)\| \geq c_\delta e^{-2 \delta m}$,
$\delta>0$, $m \geq C$.

\end{thm}

\begin{pf}

Otherwise, by (\ref {context}),(\ref{ais}),(\ref{hc}),
$|u(x+i\alpha)| \leq e^{-\delta m}$ for some $x \in \R$ and
$0 \leq i \leq m$.  Then, by
Theorem \ref{trig pol}, $\|u\|_0\leq C q_k^{Cr}e^{-\delta m}
\leq C_\delta e^{-\delta m/2}$.  This contradicts $\int u(x)=1$.
\end{pf}

Let $B(x)$ be the matrix
\be \label {Bx}
\left (\bm e^{2 \pi i \theta} u(x)&-\frac {1} {\|U(x)\|^2} \overline
{u(x-\alpha)}\\u(x-\alpha)& \frac {1} {\|U(x)\|^2}
e^{-2 \pi i \theta} \overline {u(x)} \em \right ).
\ee
By Theorem \ref {lower bound ux} and the trivial estimate $\|U\|_0
\leq C m$ we have
\be \label {2del}
\|B\|_0 \leq C_\delta e^{2 \delta m}, \quad \delta>0.
\ee
Thus
\be \la{Ber}
B(x+\alpha)^{-1} A(x) B(x)=\left (\bm e^{2 \pi i \theta}&0\\0&e^{-2 \pi i
\theta} \em \right )+\left (\bm \beta_1(x)&b(x)\\\beta_3(x)&\beta_4(x)
\em \right ),
\ee
with
\be \label {betab}
\left \|\left (\bm \beta_1(x)&b(x)\\\beta_3(x)&\beta_4(x)
\em \right ) \right \|_0 \leq C_\delta e^{4\delta m}.
\ee
Since the first column of $B$ satisfies (\ref {ais}) and $B \in
\SL(2,\C)$, (\ref {hc}) and (\ref {betab}) for appropriate $\delta$ give
\be \label {cm}
\|\beta_1\|_0,\|\beta_3\|_0,\|\beta_4\|_0 \leq
C e^{-c m}.
\ee

Taking $\Phi(x)$ the product of $B(x)^{-1}$
and a constant diagonal matrix, $\Phi(x)=D B(x)^{-1}$, where
$D=\left (\bm d&0\\0&d^{-1} \em \right )$, with $d^2=
\max \{\|\beta_3\|^{1/2}_0,e^{-m}\}$ we get
\be
\Phi(x+\alpha) A(x) \Phi(x)^{-1}=
\left (\bm e^{2 \pi i \theta}&0\\0&e^{-2 \pi i   
\theta} \em \right )+Q(x),
\ee
where, by (\ref {cm}),
(\ref {2del}),(\ref {betab}),and \eqref{d} with appropriate $\delta$,
$\|Q\|_0 \leq C e^{-c m}$ and $\|\Phi\|_0 \leq C e^m$.
This implies
\be \label {cecn}
\|A_s\|_0 \leq 2 \|\Phi\|^2_0 \leq
C e^{C m}, \quad 0 \leq s \leq c^{(1)} e^{c^{(1)} m}.
\ee

\noindent {\it Proof of Theorem \ref {growth polynomial}.}
Let $c^{(1)}$ be as in (\ref {cecn}).  For fixed $s$, let $m$ be minimal
such that $s \leq c^{(1)} e^{c^{(1)} m}$, $4 C_0 (1+|n_j|)<m<C_0^{-1}
|n_{j+1}|$, and $m=r q_k-1<q_{k+1}$
for some $j,k,r$.  By Lemma \ref {following resonance}, $m
\leq C+C \ln (1+s)$.  By (\ref {cecn}), $\|A_s\|_0 \leq C (1+s)^C$.
\qed

\subsection{Improved estimate on almost invariant sections}

In the previous argument, polynomial growth was obtained from an estimate on
a suitably chosen ``approximate Bloch wave''.  The first application of
polynomial growth is that it allows us to deal, through a bootstrap argument, with much finer approximate
Bloch waves.

\begin{lemma} \label {improved u estimate}

Let $n=|n_j|+1<\infty$ and let $N=|n_{j+1}|$.
Let $u(x)=u^I(x)$ for $I=[-C_0^{-1} N+1,C_0^{-1} N-1]$.  Define as before
$U(x)=\left (\bm e^{2\pi i \theta} u(x)\\u(x-\alpha) \em \right )$.  Then we
have the estimates
\be \label {improved u error}
A(x) \cdot U(x)=e^{2 \pi i \theta} U(x+\alpha)+
\left (\bm h(x)\\0 \em \right ),\quad \text {with} \quad
\|h\|_c \leq e^{-c N},
\ee
\be \label {improved u minoration}
\inf_{|\Im z|<c n^{-C}} \|U(x)\| \geq c n^{-C}, \quad N \geq C.
\ee
\be\label{up}
\|U\|_0\leq Cn
\ee
\end{lemma}

\begin{pf}

Estimate (\ref {improved u error}) is obtained exactly as (\ref{hc}).
For (\ref {improved u minoration}),
let $4 C_0 n<m<C_0^{-1} N$ be of the form $m=q_k-1$,
and let $J=[-[m/2],m-[m/2]]$.
Define $U^J(x)=\left (\bm e^{2\pi i \theta}
u^J(x)\\u^J(x-\alpha) \em \right )$.
Theorem \ref {almost localized} implies
\be \label {UJ1}
\|U-U^J\|_c \leq C e^{-c m}.
\ee

Arguing as in Theorem \ref {lower bound ux} we get
\be \label {UJ2}
\inf_{x \in \R/\Z} \|U^J(x)\| \geq c m^{-C}, \quad m \geq C
\ee
(the estimate is better since $|J|=q_k-1$,
and because we can use Theorem \ref {growth polynomial} instead of
(\ref{context})).
Obvious bounds on the derivative of $U^J$ then give
$\|\frac d{dz}U^J(z)\|_{c m^{-1}} \leq C m^2,$ which together with
(\ref {UJ1}) and (\ref {UJ2}) imply
\be \label {C1}
\inf_{|\Im z|<c m^{-C}} \|U(z)\| \geq c m^{-C}, \quad m \geq C^{(1)}.
\ee

By Lemma \ref {following resonance} and the Diophantine condition, if
$m>\max \{4 C_0 n,C^{(1)}\}$ is minimal of the form $m=q_k-1$ then
$m<C n^C<C_0^{-1} N$, $N \geq C$.
Together with (\ref {C1}), this gives (\ref {improved u minoration}).

Estimate \eqref{up} is immediate by Theorem  \ref {almost localized}.
\end{pf}

\section{(Complex) almost triangularization: proof of Theorem \ref
{triangula}} \label {triangulariza}

Let $U(x)$ be as in Lemma \ref {improved u estimate}.
Let $B(x)$ be defined by (\ref {Bx}).  By Lemma \ref{improved u estimate} we get the bound
\be
\|B\|_{c n^{-C}} \leq C n^C.
\ee
Here the complex extension of $B$ is the holomorphic one (not the one
given by (\ref {Bx})).  By the same computation as in (\ref{Ber}) and using Lemma \ref{improved u estimate} we get the improved estimate
\be
B(x+\alpha)^{-1} A(x) B(x)=
\left (\bm e^{2 \pi i \theta}&0\\0&e^{-2 \pi i
\theta} \em \right )+\left (\bm \beta_1(x)&b(x)\\
\beta_3(x)&\beta_4(x)\em \right )
\ee
where $\|b\|_{c n^{-C}} \leq C n^C$ and $\|\beta_1\|_{c
n^{-C}},\|\beta_3\|_{c n^{-C}},\|\beta_4\|_{c n^{-C}} \leq C e^{-c N}$.
Thus the Fourier coefficients of $b$ satisfy the
estimate
\be \la{fb}
|\hat b_k| \leq C n^C e^{-c n^{-C} |k|}.
\ee

Notice that the estimates so far cover the case of the first resonance
($n=1$), so we will assume from now on that $n>1$.

In the estimates to follow, we shift one scale back.  So take
$\tn=1+|n_{j-1}|$, $\tN=|n_j|$, and the remaining notation is clear.
Split $\tb(x)=\tb^l(x)+\tb^h(x)$ where the first term has the (low
frequency) Fourier coefficients with $|k|<\tN$ and $\tb^h$ is the rest.
We solve exactly
\be
\left (\bm 1&-\tphi(x+\alpha)\\0&1 \em \right )
\left (\bm e^{2 \pi i \theta}&\tb^l(x)\\0&e^{-2 \pi i
\theta} \em \right )
\left (\bm 1&\tphi(x)\\0&1 \em \right )=
\left (\bm e^{2 \pi i \theta}&0\\0&e^{-2 \pi i
\theta} \em \right ).
\ee
This corresponds to solving
\be
\tb^l(x)-e^{-2 \pi i \theta} \tphi(x+\alpha)+e^{2 \pi i
\theta} \tphi(x)=0,
\ee
or in terms of Fourier coefficients,
\be \la{coh}
\hat {\tphi}_k=-\hat {\tb}_k \frac {e^{-2 \pi i \theta}}
{1-e^{-2\pi i(2 \theta-k\alpha)}}
\ee
for $|k|<\tN$ and $\hat {\tphi}_k=0$ for $|k| \geq \tN$.
By (\ref{fb}), and since $\alpha\in \DC(\kappa,\tau)$
and $|k|<\tN,$ we have that $\|\tphi(x)\|_{c \tn^{-C}} \leq
C \tN^C$.  We conclude that, with $\tPhi(x)=
\left (\bm 1&-\tphi(x)\\0&1 \em \right ) \tB(x)^{-1}$ we have
\be
\tPhi(x+\alpha)A(x)\tPhi(x)^{-1}=
\left (\bm e^{2 \pi i \theta}&0\\
0&e^{-2 \pi i \theta} \em \right )+\tPsi(x)
\ee
where 
\be
\tPsi(x)=\left (\bm 1&-\tphi(x+\alpha)\\0&1 \em \right )
\left (\bm \tbeta_1(x)&\tb^h(x)\\\tbeta_3(x)&\tbeta_4(x) \em \right )
\left (\bm 1&\tphi(x)\\0&1 \em \right )
\ee
and $\|\tPhi\|_{c \tn^{-C}} \leq C \tN^C$.
Since, by (\ref{fb}), $\|\tb^h\|_{c \tn^{-C}} \leq C
\tn^C e^{-c \tn^{-C} \tN} \leq
C e^{-c \tn^{-C} \tN}$, we also have
$\|\tPsi\|_{c \tn^{-C}} \leq C e^{-c \tn^{-C} \tN}$.
It follows that
\be
\|A_s\|_{c \tn^{-C}} \leq C \tN^C, \quad 0 \leq s \leq c e^{c \tn^{-C} \tN}.
\ee

Back to the original scale, we get using Lemma \ref {following resonance}
\be \label {scale tilde n}
\|A_s\|_{c(\ln n)^{-C}} \leq C n^C, \quad 0 \leq s \leq c e^{c n (\ln n)^{-C}}.
\ee

Let $\Delta \geq n$ be minimal so that $|\hat b_k| \leq \Delta
e^{-\Delta^{-1} |k|}$, $k \in \Z$.
Let us consider a different decomposition
$b=b^r+b^l+b^h$ so that now $b^r$ has only the (resonant)
$n_j$-th Fourier coefficient,
$b^l$ has the (low frequency) Fourier coefficients with $|k| \leq \Delta^3$
(except for $n_j$) and $b^h$ has the (high frequency)
Fourier coefficients with $|k|>\Delta^3$.
By (\ref {fb}), $\Delta \leq C n^C$, and the definition of $\Delta$ gives
\be \label {bh}
\|b^h\|_{c n^{-C}} \leq C e^{-n}.
\ee

We can solve
\be
\left (\bm 1&-\phi(x+\alpha)\\0&1 \em \right )
\left (\bm e^{2 \pi i \theta}&b^l(x)\\0&e^{-2 \pi i
\theta} \em \right )
\left (\bm 1&\phi(x)\\0&1 \em \right )=
\left (\bm e^{2 \pi i \theta}&0\\0&e^{-2 \pi i
\theta} \em \right ).
\ee
As in (\ref{coh}) it follows that $\|\phi(x)\|_{c n^{-C}} \leq
C n^C$.
We conclude that, with $\Phi(x)=
\left (\bm 1&-\phi(x)\\0&1 \em \right ) B(x)^{-1}$ we obtain
\be \label {Phi}
\Phi(x+\alpha)A(x)\Phi(x)^{-1}=
\left (\bm e^{2 \pi i \theta}&0\\
0&e^{-2 \pi i \theta} \em \right )+
\left (\bm 0&b^r(x)+b^h(x)\\0&0 \em \right )+\Psi(x)
\ee
with
\be \label {phi psi}
\|\Phi\|_{c n^{-C}} \leq C n^C \quad \text {and}
\quad \|\Psi\|_{c n^{-C}} \leq C e^{-c N}.
\ee

We now estimate $b^r(x)=\hat b_{n_j} e^{2 \pi i n_j x}$.  Let
\be
W(x)=\left (\bm e^{2 \pi i \theta}&b^r(x)\\
0&e^{-2 \pi i \theta} \em \right ).
\ee
We can compute exactly
\be
W_s=\left (\bm e^{2 \pi i s\theta}&b^s(x)\\
0&e^{-2 \pi i s\theta} \em \right ).
\ee
with $|b^s(x)|=|\hat b_{n_j}\sum_{k=0}^{s-1}e^{-2\pi ik(2\theta-n_j\alpha)}|=|\hat b_{n_j}\frac{\sin \pi s
(2\theta-n_j\alpha)}{\sin \pi(2\theta-n_j\alpha)}|,\;$ if $\sin \pi (2
\theta-n_j \alpha) \neq 0$ and $|b^s(x)|=s |\hat b_{n_j}|$ otherwise.
Therefore we have
\be \label {Wsl}
\|W_s\|_0 \geq \frac {s |\hat b_{n_j}|} {10}, \quad
0 \leq s \leq \|2 \theta-n_j \alpha\|_{\R/\Z}^{-1}/10.
\ee
On the other hand, by (\ref{fb}),
\be \label {Wsu}
\|W_s\|_0 \leq 1+s |\hat b_{n_j}|\leq C (1+s) n^C, \quad s \geq 0.
\ee

Using (\ref{Phi}), we get
$A=\Phi(x+\alpha)^{-1}(W(x)+Z(x))\Phi(x)$, so that
\be \label {As}
\|A_s\|_0\geq \|\Phi\|_0^{-2} \left (\|W_s\|_0-\sum_{k=1}^s
\binom {s} {k}
\|Z\|^k_0 \max_{0 \leq j < s} \|W_j\|^{1+k}_0 \right ).
\ee
By (\ref {phi psi}) and (\ref {bh}) we have
$\|Z\|_0 \leq C e^{-n}.$ Thus, (\ref {Wsu}) and (\ref {As}) imply
\be \label {AsWs}
\|A_s\|_0 \geq cn^{-C}(\|W_s\|_0-C e^{-c n}), \quad 0 \leq s \leq c e^{c n}.
\ee

By (\ref {scale tilde n}) and (\ref {AsWs}),
$\|W_s\|_0 \leq C n^C$, $0 \leq s \leq c e^{c n (\ln n)^{-C}}$.  Together
with Lemma \ref {following resonance} and (\ref {Wsl}), we get the estimate
\be \label {br}
|\hat b_{n_j}| \leq C n^C e^{-c n (\ln n)^{-C}} \leq C e^{-c n (\ln
n)^{-C}}.
\ee
The result follows from (\ref {Phi}), (\ref {phi psi}), (\ref {bh}) and
(\ref {br}).\qed

\section{(Real) almost conjugacy to rotations: proof of Theorem \ref {di}}
\label {real conjugacies}

Let $U(x)$ be as in Lemma \ref {improved u estimate}.
Let $B(x)$ be the matrix with columns
$U(x)$ and $\overline {U(x)}$ on $\R/\Z$.  Let
$L^{-1}=\|2 \theta-n_j \alpha\|_{\R/\Z}$.  As standing hypothesis below, we
assume 
\be\label{lc}
0<L^{-1}<c.
\ee
\begin{lemma}\la{det}

For any $\eps>0,$ we have
\be \label {detmino}
\inf_{x \in \R/\Z} |\det B(x)| \geq c L^{-C}.
\ee

\end{lemma}

\begin{pf}

Recall that for any $2\times 2$ complex matrix $M$ 
with columns $V$ and $W,$

\be \la{det1}
|\det M|=\|V\|\min_{\lambda\in\C}\|W-\lambda V\|.
\ee
Notice that the minimizing $\lambda$ satisfies $\|\lambda V\|
\leq \|W\|$.

Minimize over
$\lambda \in \C$, $x \in \R/\Z$ the quantity $\|e^{-\pi i n_j x}
\overline {U(x)}-\lambda e^{\pi i n_jx} U(x)\|$.  This gives some
$\lambda_0$, $x_0$.  If the estimate does not hold then, using Lemma
\ref{improved u estimate}, Theorem \ref{growth polynomial}, (\ref{det1}), and 
Lemma \ref{following resonance}, we would have
\be \label {2piij}
\|e^{-2 \pi i s \theta} e^{-\pi i n_j x_0}
\overline {U(x_0+s\alpha)}-e^{2 \pi i s \theta} e^{\pi i n_j x_0}
\lambda_0 U
(x_0+s\alpha)\| \leq L^{-1}, \quad 0 \leq s \leq L.
\ee
This implies that
$\|e^{-\pi i n_j (x_0+s\alpha)}
\overline {U(x_0+s\alpha)}-\lambda_0 e^{\pi i n_j
(x_0+s\alpha)} U(x_0+s\alpha)\|
\leq L^{-1}+CsL^{-1} \|U\|_0 \leq L^{-1/3}$ for $0 \leq s \leq L^{1/3}$
(recalling \eqref{up}).

As in the proof of Lemma \ref{improved u estimate}, let
$m>(\ln L)^2>4 C_0 n$ (here the second inequality is a
  consequence of \eqref{lc}) be minimal
of the form $m=q_k-1$, let
$J=[-[m/2],m-[m/2]]$ and let $U^J(x)=\left (\bm e^{2\pi i \theta}
u^J(x)\\u^J(x-\alpha) \em \right )$.  By Lemma \ref
{following resonance} and the Diophantine condition,
$m<C (\ln L)^C<C_0^{-1} N$.
Theorem \ref {almost localized} then gives
$\|U-U^J\|_c \leq L^{-1}$.
Using Theorem \ref{trig pol} and Theorem \ref{growth polynomial}, we obtain 
$\|e^{-\pi i n_j x}
\overline {U^J(x)}-e^{\pi i n_j x} \lambda_0 U^J(x)\|_0 \leq C (\ln L)^C
(L^{-1/3}+2 L^{-1})$.
Thus
\be \label {1/4}
\sup_{x \in \R/\Z} \|e^{-\pi i n_j x}
\overline {U(x)}-e^{\pi i n_j x} \lambda_0 U(x)\| \leq L^{-1/4}.
\ee
Substituting $x_1=x_0+s\alpha$ in (\ref {2piij}) and taking $s=[L/2]$,
we get $\|i e^{-\pi i n_j x_1} \overline
{U(x_1)}+i e^{\pi i n_j x_1} \lambda_0 U(x_1)\| \leq L^{-1}+C L^{-1}
\|U\|_0$, so that, by (\ref {1/4}) and
\eqref{up}, $\|U(x_1)\| \leq L^{-1/5}$.

However, (\ref{improved u minoration}) and Lemma \ref {following
  resonance} imply
$\|U(x_1)\| \geq c (\ln L)^{-C}$, giving a contradiction.

\end{pf}

\begin{lemma}

Let $x_0 \in \R/\Z$.  Then
\be \la{detvar}
\sup_{|\Im z|<c} |\det B(z)-\det B(x_0)| \leq C e^{-c N}.
\ee

\end{lemma}

\begin{pf}

By (\ref{improved u error}), $A(x) 
B(x)=\left (\bm e^{2 \pi i \theta}&0\\0&e^{-2 \pi i \theta} \em \right )B(x+\alpha)+
\left (\bm h(x)\\0 \em \right ),\quad \|h\|_c \leq e^{-c N},$
thus $|\det B(x_0+\alpha)-\det B(x_0)| \leq C e^{-cN}.$ This gives
$|\det B(x_0+k\alpha)-\det B(x_0)| \leq C e^{-c N}$,
$0 \leq k \leq 4 N$.  The
function $x \mapsto \det B(x_0+k\alpha)-\det B(x)$ is a trigonometrical
polynomial of essential degree bounded by $4 N$.  By Theorem \ref {trig
pol}, $|\det B(x)-\det B(x_0)| \leq C e^{-c N}$, $x \in \R/\Z$.
By Theorem \ref{almost localized}, $\|U(z)\|_c
\leq C e^{C n}$,\footnote {Actually, even a trivial bound like
$\|U(z)\|_c \leq C e^{C N}$ would do here.} so 
$z \mapsto \det B(z)-\det B(x_0)$
is bounded by $C e^{C n}$ over
$|\Im z|<c^{(2)}$.  By the Hadamard three-circle theorem, 
\be
\ln \sup_{|\Im z|=\delta c^{(2)}} |\det B(z)-\det B(x_0)| \leq (1-\delta)
\ln \sup_{|\Im
  z|=0} +\delta\ln \sup_{|\Im
  z|=c^{(2)}} \leq -cN, \quad 0 \leq \delta<c.
\ee

Thus $|\det B(z)-\det B(x_0)| \leq C e^{-c N}$, $|\Im z|<c$.
\end{pf}

\begin{thm} \la{det2}

We have
\be
\inf_{|\Im z|<c} |\Im \det B(z)| \geq c L^{-C}\geq \frac C{N^C}.
\ee

\end{thm}

\begin{pf}

Notice that $\Re \det B(x)=0$, $x \in \R/\Z$.  The result now
follows from (\ref {detmino}) and (\ref {detvar})
(taking into account Lemma \ref {following resonance}).
\end{pf}

\begin{rem}

Optimizing the method (and using other estimates obtained in this paper)
gives $|\det B(z)| \geq c_\delta L^{-1-\delta}$, $|\Im z|<c$.
This will be explored and used in \cite{ajlp}.

\end{rem}

Take now $S=\Re U$, $T=\Im U$ on $\R/\Z$.
Then $B=\left (\bm S & \pm T \em
\right )\left (\bm 1&1\\\pm i &\mp i \em \right ).$  Let $ W_1$
be the matrix with columns $S$ and $\pm T$, so to have $\det W_1>0$.
Since $\left (\bm 1&1\\\pm i &\mp i \em \right )\left (\bm e^{2 \pi
  i \theta}&0\\0&e^{-2 \pi i \theta} \em \right )\left (\bm 1&1\\\pm
i &\mp i \em \right )^{-1}=R_{\mp\theta},$ by (\ref{improved u error})
we have
\be
\|A(x) \cdot  W_1(x)-
W_1(x+\alpha) \cdot R_{\mp \theta}\| \leq C e^{-c N}, \quad |\Im x|<c
\ee
(the complex extension considered here is the holomorphic one).
Let $W(x)=|\det B(x)/2|^{-1/2} W_1(x)$ on $\R/\Z$,
so to have $\det W=1$ (by Theorem \ref {det2}, there is no problem with
branching when extending $|\det B(x)|^{-1/2}$ to $|\Im x|<c$).  Then
\be
\left \|A(x) \cdot W(x)-\frac {|\det B(x+\alpha)|^{1/2}} {|\det B(x)|^{1/2}}
W(x+\alpha) \cdot R_{\mp \theta} \right \| \leq C e^{-c N},
\quad |\Im x|<c.
\ee

A combination of (\ref{detvar}) and Theorem \ref{det2} gives 
\be \la{fdet}
\left |\frac {|\det B(x+\alpha)|^{1/2}} {|\det B(x)|^{1/2}}-1 \right | \leq
C e^{-c N}, \quad |\Im x|<c.
\ee
Also, using Theorems \ref{almost localized} and \ref{det2} we get
$\|W\|_c \leq CN^Ce^{cn}.$ Thus
\be
\|A(x) \cdot W(x)-W(x+\alpha) \cdot R_{\mp \theta}\| \leq
Ce^{-c N}, \quad |\Im x|<c.
\ee

To conclude, we need to show that $|\deg W| \leq C n$.
Set $\tilde U(x)= e^{\pi i n_jx} U(x),$ and $\tilde S=\Re \tilde U,\;
\tilde T=\Im \tilde U,\;
\tilde W_1 = \left (\bm \tilde S & \pm \tilde T \em
\right ),\; \tilde W(x)=|\det B(x)/2|^{1/2} \tilde W_1(x)$ on $\R/2\Z$.
Then $\tilde W(x)=W(x)R_{\mp\frac
  {n_j x}2},$ so $\deg \tilde W-2\deg W=\mp n_j$, where $\deg \tilde W$ is
the degree of $\tilde W:\R/2\Z \to \SL(2,\R)$ and $\deg W$ is
understood in the usual sense, as the degree of $W:\R/\Z \to \SL(2,\R)$.

The degree of $\tilde W:\R/2\Z \to \SL(2,\R)$ is the same as the degree
of any of its columns, considered as maps $\R/2\Z \to \R^2 \setminus \{0\}$.
Thus we need to show that the degree
of $M:\R/2\Z \to \R^2 \setminus \{0\}$,
for either $M=\tilde S$ or for $M=\tilde T$ satisfies $|\deg M| \leq Cn$.

Notice that 
\be \label{st} \|\int_{\R/2\Z} e^{-\pi in_j x}(\tilde S(x)+i \tilde T(x))
dx\|=\|\int_{\R/2\Z} U(x) dx\| \geq 2,\ee since $\hat u_0=1.$ Select
$M(x)=\tilde S(x)$ or $M(x)=\tilde T(x)$ so that $\int_{\R/2\Z} \|M(x)\| \geq 1$.
Then, by (\ref{improved u error}), for $|\Im x|<c $ we have
$\|A(x) \cdot \tilde U(x)-
e^{\pi i (2\theta-n_j\alpha)} \tilde U(x+\alpha)\|_0 \leq C e^{-c N}$. 
By Lemma \ref {following resonance} and using \eqref{up}
and $L^{-1}\leq e^{-\epsilon_0|n|},$
we have $|\|A(x) \cdot M(x)\|-\|M(x+\alpha)\|| \leq CL^{-1}n\leq C L^{-c}$,
$x \in \R/2\Z$.  Using Theorem
\ref{growth polynomial}, we obtain, arguing as in the proof of Lemma
\ref{improved u estimate} 
\be \la{infm}
\inf_{|\Im x|<cn^{-C}} \|M(x)\| \geq c n^{-C}.
\ee

This already gives
the bound $|\deg M| \leq C n^C$ by the obvious derivative
estimate.  To obtain a linear bound, let
$\tilde M(x)$ be a vector obtained
by cutting off the Fourier modes of $M(x)$ with $|\frac {k} {2}|>\Delta$,
where $\Delta$ is chosen minimal so that
$\|\tilde M(x)-M(x)\|<\|M(x)\|$, $x \in \R/2\Z$.
Then by Theorem \ref {almost localized} and (\ref {infm}),
$\Delta<C n$.  By Rouche's Theorem, the degree of $M$ is
the same as the degree of $\tilde M$.
Consider now a coordinate
of $x \mapsto \tilde M(2x)$ which is not identically vanishing.  It is a
trigonometric polynomial of essential degree at most $C n$, and since it not
identically vanishing, it has at most $C n$ zeroes in
$\R/\Z$.  It follows that $|\deg \tilde M|$ (and thus $|\deg M|$) is bounded
by $C n$, and we conclude that $|\deg W| \leq Cn$.\qed

\appendix

\section{The perturbative theory of Eliasson}

The following result is due to Hakan Eliasson.

\begin{thm}\la{eli}

Let $\alpha \in \DC(\kappa,\tau)$, and let $A:\R/\Z \to \SL(2,\R)$ be
analytic.  Assume that
$\|A-A_*\|_\epsilon<C^{-1}_0 \epsilon^{r_0}$ where
$A_* \in \SL(2,\R)$, $0<\epsilon<1$, $r_0=r_0(\tau)$ and
$C_0=C_0(\kappa,\tau,\|A_*\|)$.  Then
\begin{enumerate}
\item If $\rho(\alpha,A)$ is either Diophantine or rational with
respect to $\alpha$ then $(\alpha,A)$ is analytically reducible,
\item If $(\alpha,A)$ is not hyperbolic then $\sup \frac {1} {k}
\|A_k\|_0 \leq C_0$,
\item If $\rho(\alpha,A)$ is not rational then
$\lim \frac {1} {k} \|A_k\|_0=0$.
\end{enumerate}

\end{thm}

Actually what is considered in \cite {E} is the case of continuous time, and
for cocycles of Schr\"odinger type (he also treats the case of several
frequencies).  The considerations for discrete time
are similar and are carried out in the thesis of his student
Sana Ben Hadj Amor\cite{Am}.  In \cite {E}, Eliasson goes on to establish absolutely
continuous spectrum for all phases for the associated Schr\"odinger
operators, while
the proof of $1/2$-H\"older continuity of the integrated density
of states is the main result of Amor's thesis.

The KAM scheme of Eliasson
has been extended to the $C^\infty$ case \cite {AK}.  This generalization
motivated the introduction of the concept of ``almost reducibility'' in the
$C^\infty$ case.

\subsection{Quick reduction to the perturbative regime}

Reduction to the perturbative regime of Eliasson is a much less subtle
result than showing almost reducibility.  It can be concluded from almost
localization just after \S \ref {growth}, since is
clearly\footnote{Recall that if $\theta$ is non-resonant then by  Theorems \ref{localreduc} and
  \ref{almost localized} $(\alpha,A)$ is conjugate to a
constant cocycle (which is obviously in Eliasson's perturbative regime).} implied by the
following result.

\begin{thm}\label{a1}

If $\theta$ is resonant, there exists arbitrarily small
$\epsilon>0$, $-1 \leq \kappa \leq 1$,
and $B:\R/\Z \to \PSL(2,\R)$ analytic such
that $\|B\|_\epsilon \leq C \epsilon^{-C}$ and such that
$\|B(x+\alpha)A(x)B(x)^{-1}-A_*\|_\epsilon \leq
C e^{-c/\epsilon^c}$, where
$A_*=\left (\bm 1&\kappa\\0&1 \em \right )$.

\end{thm}

\comm{
Recall that if $\theta$ is non-resonant then $(\alpha,A)$ is conjugate to a
constant cocycle (which is obviously in Eliasson's perturbative regime).
We will now show that if $\theta$ is non-resonant then we
can conjugate it to a cocycle in Eliasson's perturbative regime.
}

\begin{pf}

The beginning of the proof coincides with a small part of the argument
used in the proof of Theorem \ref{di} but we briefly repeat it here for the
readers' convenience. Let $n=|n_j|+1$, $N=|n_{j+1}|$,
and let $I=[-C_0^{-1} N+1,C_0^{-1} N-1]$.
Let $u=u^I$, and define $U(x)=\left (\bm e^{2 \pi i \theta}
u(x)\\u(x-\alpha) \em \right )$ as usual. Let 
$\tilde U(x)=e^{\pi i n_j x} U(x)$ and let $\tilde S(x)=\Re \tilde U(x)$,
$\tilde T(x)=\Im \tilde U(x)$.  Since (\ref{st}) holds, we can choose
$M=\tilde S$ or $M=\tilde T$ so that $\int_{\R/2\Z} \|M(x)\| \geq
1$.  Notice that by Lemma \ref{improved u estimate} and the definition
of resonance \be\label{ai}A(x) \cdot M(x)=M(x+\alpha)+O(e^{-cn}),\;\;|\Im x|<c.\ee
Then, arguing as in Lemma \ref {improved u estimate}, we have \eqref{infm}.

Let $W$ be the matrix with columns $M$ and $\frac {1} {\|M\|^2} R_{1/4} M$. 
Then by \eqref{ai} and \eqref{infm} we have $W(x+\alpha)^{-1}A(x)W(x)=\left (\bm 1&0\\0&1\em \right )+\left (\bm
\beta_1(x)&\beta_2(x) \\ \beta_3(x)&\beta_4(x) \em \right )$ with
$\|\beta_1\|_{c n^{-C}},\|\beta_3\|_{c n^{-C}}, \|\beta_4\|_{c n^{-C}} \leq
C e^{-cn}$ and $\|\beta_2\|_{c n^{-C}} \leq C n^C$.  Solve the cohomological
equation $\phi(x+\alpha)+\beta_2(x)-\phi(x)=\int_0^2 \beta_2(x) dx/2=b$, with
$\int_0^2 \phi(x) dx=0$.  Then, by writing out the Fourier
coefficients and since $\alpha\in\DC$ we get $\|\phi\|_{c n^{-C}} \leq C n^C$.  Let
$\Phi(x)=\left (\bm 1&\phi(x)\\0&1 \em \right ) W(x)^{-1}$, then
$\Phi(x+\alpha)A\Phi (x)^{-1}=\left (\bm 1&b\\0&1 \em \right )+H$
where $\|H\|_{cn^{-C}}\leq Ce^{-cn}.$ 
Let $d^2=\min \{|b|^{-1},1\}$ and let $D=\left (\bm
d&0\\0&d^{-1} \em \right )$.  Then $B(x)=D \Phi(x+\alpha)
A(x) \Phi(x)^{-1} D^{-1}$ is such that
\be
\|B-\left (\bm 1&\kappa\\0&1 \em \right )\|_{c n^{-C}} \leq Ce^{-cn},
\ee
for some $-1 \leq \kappa \leq 1$.
\end{pf}


\begin{thebibliography}{BKNS}

\bibitem[Am]{Am} Amor, S.
Op\'erateurs de Schr\"odinger quasi-p\'eriodiques uni-dimensionnels. 
Thesis, Universit\'e Paris 7, 2006.

\bibitem[AA]{AA} Aubry, S.; Andr\'e, G.
Analyticity breaking and
Anderson localization in incommensurate lattices.  Group theoretical methods
in physics (Proc. Eighth Internat. Colloq., Kiryat Anavim, 1979),  pp.
133--164, Ann. Israel Phys. Soc., 3, Hilger, Bristol, 1980.

\bibitem[A]{a} Avila, A. Absolutely continuous spectrum for the almost
  Mathieu operator with subcritical coupling. In preparation.

\bibitem[AD]{ad} Avila, A.;Damanik, D. Absolute continuity of the
integrated density of states for the almost Mathieu operator with
non-critical coupling.  Inventiones Mathematicae 172 (2008), 439-453.

\bibitem[AJ1]{AJ} Avila, A.;Jitomirskaya, S.
The Ten Martini Problem.  Preprint (www.arXiv.org).  To appear in Annals of Mathematics.

\bibitem[AJ2]{ajlp}  Avila, A.;Jitomirskaya, S. In preparation.

\bibitem[AK1]{ak0} Avila, A.; Krikorian, R. Reducibility or non-uniform hyperbolicity for quasiperiodic Schr\"odinger
cocycles.  Annals of Mathematics 164 (2006), 911-940.

\bibitem[AK2]{AK} Avila, A.; Krikorian, R.
Some local and semi-local results for quasiperiodic Schr\"odinger cocycles.
In preparation.




\bibitem[BLT]{blt}Bellissard, J.; Lima, R.; Testard, D. Almost periodic Schr\"odinger operators. Mathematics + physics. Vol. 1, 1--64, World Sci. Publishing, Singapore, 1985
\bibitem[Be]{ber} Berezanskii, Y.
Expansions in eigenfunctions of selfadjoint operators. Transl. Math. Monogr.,
Vol. 17. Providence, RI: Am. Math. Soc. 1968.



\bibitem[B1]{B1}  Bourgain, J. H\"older regularity of
   integrated density of states for the almost Mathieu operator in a
   perturbative regime. Lett. Math. Phys. 51 (2000), no. 2, 83--118.

\bibitem[B2]{B3} Bourgain, J. On the spectrum of lattice Schr\"odinger
operators with deterministic potential. II. Dedicated to the memory of Tom
Wolff.  J. Anal. Math. 88 (2002), 221-254.

\bibitem[B3]{B2} Bourgain, J. Green's function
estimates for lattice Schr\"odinger operators and applications. Annals
of Mathematics Studies, 158. Princeton University Press, Princeton,
NJ, 2005. x+173 pp.






\bibitem[AS]{as} Avron, Joseph; Simon, Barry Almost periodic Schr�inger
operators. II. The integrated density of states.  Duke Math. J.  50  (1983),
no. 1, 369--391.





\bibitem[BJ1]{BJ1} Bourgain, J.; Jitomirskaya, S.
Absolutely continuous spectrum for 1D quasiperiodic operators.
Invent. Math. 148 (2002), no. 3, 453--463.

\bibitem[BJ2]{BJ2} Bourgain, J.; Jitomirskaya, S. Continuity of the Lyapunov
exponent for quasiperiodic operators with analytic potential. Dedicated to
David Ruelle and Yasha Sinai on the occasion of their 65th birthdays.  J.
Statist. Phys.  108  (2002),  no. 5-6, 1203--1218.




\bibitem[CL]{Cy} Carmona, R.,  Lacroix, J.  Spectral Theory of Random
Schr\"odinger Operators.  Boston, MA:  Birkhauser 1990.
\bibitem[CEY]{CEY} Choi, Man Duen; Elliott, George A.; Yui, Noriko Gauss
polynomials and the rotation algebra.  Invent. Math.  99  (1990),  no. 2,
225--246.










\bibitem[DJ]{jo}De Concini, Corrado; Johnson, Russell A. The algebraic-geometric AKNS potentials. Ergodic Theory Dynam. Systems 7 (1987), no. 1, 1--24.

\bibitem[DeS]{DS} Deift, P.; Simon, B.
Almost periodic Schr\"odinger operators. III. The absolutely continuous
spectrum in one dimension.  Comm. Math. Phys.  90  (1983),  no. 3, 389--411.

\bibitem[DiS]{ds}Dinaburg, E.; Sinai, Ya. The one-dimensional 
Schr\"{o}dinger equation with a quasi-periodic potential. Funct. Anal.
Appl.{\bf 9} (1975), 279--289.
\bibitem[E]{E} Eliasson, L. H.
Floquet solutions for the $1$-dimensional quasi-periodic
Schr\"odinger equation.
Comm. Math. Phys. 146 (1992), no. 3, 447--482.

\bibitem[GS]{GS} Goldstein, M.; Schlag, W.
Fine properties of the integrated density of states and a quantitative
separation property of the Dirichlet eigenvalues.  Preprint (www.arXiv.org).
\bibitem[GJLS]{gjls}Gordon, A. Y.; Jitomirskaya, S.; Last, Y.; Simon,
  B. Duality and singular continuous spectrum in the almost Mathieu equation. Acta Math. 178 (1997), no. 2, 169--183.




\bibitem[H]{H} Herman, Michael-R. Une m\'ethode pour minorer les
exposants de Lyapounov et quelques exemples montrant le caract\`ere local
d'un th\'eor\`eme
d'Arnol'd et de Moser sur le tore de dimension $2$.
Comment. Math. Helv.  58  (1983),  no. 3, 453--502.









\bibitem[J]{J} Jitomirskaya, Svetlana Ya. Metal-insulator transition for the
almost Mathieu operator.  Ann. of Math. (2)  150  (1999),  no. 3,
1159--1175.

\bibitem[JL]{JL} Jitomirskaya, Svetlana Ya.; Last, Yoram
   Anderson localization for the almost Mathieu equation. III.
   Semi-uniform localization, continuity of gaps, and measure of the
   spectrum. Comm. Math. Phys. 195 (1998), no. 1, 1--14.









\bibitem[JS]{js}Jitomirskaya, S.; Simon, B. Operators with singular continuous spectrum. III. Almost periodic Schr\"odinger operators. Comm. Math. Phys. 165 (1994), no. 1, 201--205.

\bibitem[JM]{JM} Johnson, R.; Moser, J. The rotation number for almost
periodic potentials.  Comm. Math. Phys.  84  (1982), no. 3, 403--438.










\bibitem[LS]{LS} Last, Yoram; Simon, Barry
Eigenfunctions, transfer matrices, and absolutely continuous spectrum of
one-dimensional Schr\"odinger operators.  Invent. Math.  135  (1999),  no. 2,
329--367.

\bibitem[MP]{hyp} Moser, J. ; P\"oschel, J.
An extension of a result by Dinaburg and Sinai on quasi-periodic potentials. 
Comment. Math. Helv. 59 (1984), 39-85.





\bibitem[P1]{P} Puig, Joaquim Cantor spectrum for the almost Mathieu
operator.  Comm. Math. Phys.  244  (2004),  no. 2, 297--309.
\bibitem[P2]{puig2} Puig, Joaquim A nonperturbative Eliasson's reducibility theorem. Nonlinearity 19 (2006), no. 2, 355--376.






\bibitem[S1]{sim} B. Simon, Schr\"{o}dinger semigroups, Bull. AMS {\bf 7} (1982),
 447--526 

\bibitem[S2]{S} Simon, Barry Schr\"odinger operators in the twenty-first
century.  Mathematical physics 2000,  283--288, Imp. Coll. Press, London,
2000.








\end{thebibliography}
\end{document}